\documentclass[10pt]{amsart}
\usepackage{mathptmx}
\usepackage{amsmath}     
\usepackage{amssymb}
\usepackage{array}
\usepackage{geometry}
\usepackage[bookmarks=true,colorlinks=true, pdfstartview=FitV, linkcolor=black, citecolor=blue, urlcolor=black]{hyperref}
\usepackage{movie15}

\usepackage{color}
\definecolor{DarkRed}{rgb}{0.55,.00,0.2}
\definecolor{DarkGrey}{rgb}{0.35,.35,0.35}

\theoremstyle{definition}

\theoremstyle{remark}

\numberwithin{equation}{section}

\newcommand{\e}{{\rm e}}

\begin{document}

\title{Orthogonal polynomials  with ultra-exponential weight functions: an explicit solution to the Ditkin-Prudnikov problem}

\author{S. Yakubovich}
\address{Department of Mathematics, Fac. Sciences of University of Porto,Rua do Campo Alegre,  687; 4169-007 Porto (Portugal)}
\email{ syakubov@fc.up.pt}
\thanks{\thanks{The work was partially supported by CMUP [UID/MAT/00144/2013], which is funded by FCT(Portugal) with national (MEC) and European structural funds through the programs FEDER, under the partnership agreement PT2020, and Project STRIDE - NORTE-01-0145-FEDER- 000033, funded by ERDF - NORTE 2020.  The author thanks Marco Martins Afonso for  necessary numerical calculations and verifications of some formulas. 
 }}

\subjclass[2000]{ 33C47,  33C45, 33C10, 44A15, 42C05 }

\date{\today}


\keywords{Orthogonal polynomials, modified Bessel functions, Meijer $G$-function, Mellin transform, associated Laguerre polynomials, multiple orthogonal polynomials}

\begin{abstract}  New sequences of orthogonal polynomials with ultra-exponential weight functions are discovered.  In particular, it gives an explicit solution to the  Ditkin-Prudnikov problem (1966).   The 3-term recurrence relations, explicit representations, generating functions and Rodrigues-type formulae  are derived. The method is based on differential properties of the involved special functions and  their representations in terms of the Mellin-Barnes and Laplace integrals. A notion of the composition polynomial orthogonality is introduced. The corresponding advantages of this orthogonality to discover new sequences of  polynomials and their relations to the corresponding multiple orthogonal polynomial ensembles are shown.   \\
\end{abstract}

\maketitle

\markboth{\rm \centerline{ S. Yakubovich}}{Orthogonal polynomials  with ultra-exponential weight functions}



\section{Introduction and preliminary results}

Throughout  the text, $\mathbb{N}$ will denote the set of all positive integers, $\mathbb{N}_{0}=\mathbb{N}\cup \{0\}$, whereas $\mathbb{R}$ and $\mathbb{C}$  the field of the real and complex numbers, respectively. The notation $\mathbb{R}_{+}$   corresponds to the set of all positive real numbers. The present investigation is primarily targeted at analysis of sequences of orthogonal polynomials with respect to the weight functions related to the modified Bessel functions  of the second kind or Macdonald functions $K_\nu(x)$ \cite{Bateman}, Vol. II. The problem was posed by Ditkin and Prudnikov in the seminal work of 1966 \cite{DitkinPrudnikov} to find a new sequence of orthogonal polynomials $\left(P_n\right)_{n \in \mathbb{N}_{0}}$, satisfying the orthogonality conditions

$$\int_{0}^{\infty} 2 K_{0}(2\sqrt{x}) P_{m} (x) P_{n}(x) dx = \delta_{n,m} \ , \quad n,m\in\mathbb{N}_{0},\eqno(1.1)$$
where $\delta_{n,m}$ represents the Kronecker symbol, and  related to the weight $2K_0(2\sqrt x)$ which can be defined in terms of the Mellin-Barnes integral (see \cite{PrudnikovMarichev}, relation (8.4.23.1), Vol. III

$$
	2K_{0}(2\sqrt{x}) = \frac{1}{2\pi i} \int_{\gamma-i\infty}^{\gamma+i\infty} \Gamma^{2}(s) x^{-s} ds\ , \quad x, \gamma \in \mathbb{R}_{+},\eqno(1.2)
$$
where $\Gamma(z)$ is the Euler gamma-function \cite{Bateman}, Vol. I. The first four  polynomials are

$$P_0(x)= 1,\quad  P_1(x)= {1\over \sqrt 3} (x-1), \quad P_2(x)= \sqrt{3\over 41} \left( {x^2\over 4}- {8\over 3} x + {5\over 3}\right),$$

$$P_3(x)= \sqrt{41\over 2841} \left( {x^3\over 36}- {177\over 164} x^2 + {267\over 41} x - {131\over 41}\right).$$
Later in 1993 \cite{PrudnikovProblem} Prudnikov formulated the problem in terms of  more general ultra-exponential weight functions $\rho_{0,k-1}, \ k \in \mathbb{N}$ (see Definition 1 below),  and in \cite{AsscheYakubov2000} it was announced in terms of the scaled Macdonald function 

$$\rho_\nu(x)= 2 x^{\nu/2} K_\nu(2\sqrt x),\ x \in  \mathbb{R}_{+},\quad \nu \ge 0.\eqno(1.3)$$
This function has  the Mellin-Barnes integral representation in the form
$$
\rho_\nu(x)=  \frac{1}{2\pi i} \int_{\gamma-i\infty}^{\gamma+i\infty} \Gamma(\nu+s) \Gamma (s) x^{-s} ds\ , \quad x, \gamma \in \mathbb{R}_{+},\eqno(1.4)
$$
and more general  ultra-exponential weight functions  can be represented, in turn,  in terms of Meijer $G$-functions \cite{YaL}.  Namely, the problem is to find a sequence of orthogonal polynomials $\left(P_n^\nu\right)_{n\in \mathbb{N}_0} \ (P_n^0 \equiv P_n)$, satisfying the following orthogonality conditions

$$\int_0^\infty  P_n^\nu(x)P_m^\nu(x) \rho_\nu(x) dx=  \delta_{n,m} , \quad n,m\in\mathbb{N}_{0}.\eqno(1.5)$$
As it was shown in \cite{AsscheYakubov2000} and \cite{Cous} it is more natural to investigate multiple orthogonal polynomials for two Macdonald weights $\rho_\nu$ and $\rho_{\nu+1}$ since it gives explicit formulas, differential properties, recurrence relations and Rodrigues formula. Nevertheless, an original problem still attracts to understand the nature of such  polynomial sequences and their relation to classical systems of orthogonal polynomials and associated multiple orthogonal polynomial ensembles.

On the other hand, the operational calculus associated to the differential operator $\frac{d}{dt}$ gives rise to the Laplace transform 

$$F(x)= \int_0^\infty e^{-xt} f(t) dt, \quad x \in  \mathbb{R}_{+},\eqno(1.6)$$
having the exponential function as a kernel, which is the weight function for classical Laguerre polynomials \cite{Chi}, being represented  in terms of the Mellin-Barnes integral \cite{PrudnikovMarichev}, relation (8.4..3.1), Vol. III
$$
	\e^{-x} = \frac{1}{2\pi i} \int_{\gamma-i\infty}^{\gamma+i\infty} \Gamma(s) x^{-s} ds	\ , \quad x, \gamma  \in \mathbb{R}_{+}.\eqno(1.7) 
$$
Meanwhile,  the operator  $\frac{d}{dt}t\frac{d}{dt}$ which is also called the Laguerre derivative \cite{Dattoli},  leads to the Meijer transform \cite{YaL},  involving the weight $2K_{0}(2\sqrt{x})$ which is given by (1.2), namely, 

$$G(x)=  \int_{0}^{\infty} 2 K_{0}(2\sqrt{xt}) g(t) dt,\quad  x \in \mathbb{R}_{+}.\eqno(1.8)$$
This transform is an important example of the so-called Mellin type convolution transforms, which are extensively investigated in \cite{YaL}.  Moreover,  in the sequel we will employ the Mellin  transform technique developed in \cite{YaL} in order to  investigate various properties of the scaled Macdonald functions and more general ultra-exponential weights.   Precisely, the Mellin transform is defined, for instance, in  $L_{\mu, p}(\mathbb{R}_+),\ 1 \le  p \le 2$ (see details in \cite{Tit}) by the integral  
$$f^*(s)= \int_0^\infty f(x) x^{s-1} dx, \quad s \in \mathbb{C},\eqno(1.9)$$
 being convergent  in mean with respect to the norm in $L_q(\mu- i\infty, \nu + i\infty),\ \mu \in \mathbb{R}, \   q=p/(p-1)$.   Moreover, the  Parseval equality holds for $f \in L_{\mu, p}(\mathbb{R}_+),\  g \in L_{1-\mu, q}(\mathbb{R}_+)$
$$\int_0^\infty f(x) g(x) dx= {1\over 2\pi i} \int_{\mu- i\infty}^{\mu+i\infty} f^*(s) g^*(1-s) ds.\eqno(1.10)$$
The inverse Mellin transform is given accordingly
 $$f(x)= {1\over 2\pi i}  \int_{\mu- i\infty}^{\mu+i\infty} f^*(s)  x^{-s} ds,\eqno(1.11)$$
where the integral converges in mean with respect to the norm  in   $L_{\mu, p}(\mathbb{R}_+)$
$$||f||_{\mu,p} = \left( \int_0^\infty  |f(x)|^p x^{\mu p-1} dx\right)^{1/p}.\eqno(1.12)$$
In particular, letting $\mu= 1/p$ we get the usual space $L_p(\mathbb{R}_+; \  dx)$.   Recalling the Meijer transform (1.8) one can treat it as an analog of the Laplace transform (1.6) in the operational calculus associated with the Laguerre derivative. Consequently, the corresponding analog of the classical Laguerre polynomials would be important to investigate, discovering the mentioned Ditkin-Prudnikov polynomial sequence. Finally,  we note in this section that  in \cite{Filipa} some non-orthogonal polynomial systems were investigated which share the same canonical regular form with Ditkin-Prudnikov polynomial sequence $\left(P_n\right)_{n\in \mathbb{N}_0}$.  An analogous relation occurs, for instance,  between the Bernoulli polynomials, which also happen to be non-orthogonal, and the (orthogonal) Legendre polynomials.

\section{Properties of the scaled Macdonald functions}

We begin with

{\bf Definition 1}. {\it Let $x,\gamma  \in \mathbb{R}_+,  \  \nu \ge 0, \ k \in \mathbb{N}_0$.  The function  $\rho_{\nu,k}(x)$ is called the ultra-exponential weight function and it is expressed in terms of the following Mellin-Barnes integral}

$$\rho_{\nu,k}(x)=  \frac{1}{2\pi i} \int_{\gamma-i\infty}^{\gamma+i\infty} \Gamma(\nu+s) \left[ \Gamma (s) \right]^k  x^{-s} ds.\eqno(2.1)$$
It is easily seen from the reciprocal formulas (1.9), (1.11) for the Mellin transform that the case $k=0$ corresponds to the weight function $\rho_{\nu,0}(x)= x^\nu e^{-x}$, which is related to the classical associated Laguerre polynomials $L_n^\nu(x)$ \cite{Chi}

$$\int_0^\infty L^\nu_n(x) L_m^\nu(x) e^{-x} x^\nu dx =  \delta_{n,m} , \quad n,m\in\mathbb{N}_{0}.\eqno(2.2)$$
and $k=1$ gives the function $\rho_{\nu,1}\equiv \rho_\nu$, which  is associated with the Prudnikov polynomials $P_n^\nu$ under orthogonality conditions (1.5).  As it was mentioned above the weights $\rho_{\nu,k}$ can be expressed in terms of the Meijer $G$-functions (cf. \cite{Arno}). Concerning the scaled Macdonald function $\rho_\nu$, we employ the Parseval equality (1.10) to the integral (1.4) to derive the Laplace integral representation for this weight function which will be used in the sequel. In fact, we obtain

$$\rho_\nu(x)= \int_0^\infty t^{\nu-1} e^{-t - x/t} dt,\quad  x >0,\ \nu \in \mathbb{R}.\eqno(2.3)$$ 
The direct Mellin transform (1.9) gives the moments of $\rho_\nu$. Precisely, we obtain

$$\int_0^\infty  \rho_\nu(x) x^\mu dx = \Gamma(\mu+\nu+1)\Gamma(\mu+1).$$
Moreover, the asymptotic behavior of the modified Bessel function at infinity and near the origin \cite{Bateman}, Vol. II gives the corresponding values for the scaled Macdonald function $\rho_\nu,\ \nu \in \mathbb{R}$.  Precisely, we have
$$\rho_\nu (x)= O\left( x^{(\nu-|\nu|)/2}\right),\  x \to 0,\ \nu\neq 0, \quad  \rho_0(x)= O( \log x),\ x \to 0,\eqno(2.4)$$

$$ \rho_\nu(x)= O\left( x^{\nu/2- 1/4} e^{- 2\sqrt x} \right),\ x \to +\infty.\eqno(2.5)$$
Returning to the Mellin-Barnes integral (1.4), we multiply both sides of this equality by $x^{-\nu}$ and then differentiate with respect to $x$ under the integral sign. This is possible  via the absolute and uniform convergence by $x \ge x_0 >0$, which can be established using the Stirling asymptotic formula for the gamma-function \cite{Bateman}, Vol. I. Therefore we deduce

$${d\over dx} \left[ x^{-\nu} \rho_\nu(x)\right] =  - \frac{1}{2\pi i} \int_{\gamma-i\infty}^{\gamma+i\infty} \Gamma(\nu+s+1) \Gamma (s) x^{-s-\nu-1} ds,$$ 
where the reduction formula $\Gamma(z+1)= z\Gamma(z)$ for the gamma-function is applied. Multiplying the latter equality by $x^{\nu+1}$ and differentiating again, we involve a simple change of variables and the analyticity on the right half-plane ${\rm Re} s > 0$ of the integrand to end up with the second order differential equation for $\rho_\nu$

$$ {d\over dx} \left[  x^{\nu+1} {d\over dx} \left[ x^{-\nu} \rho_\nu(x)\right] \right] = \rho_\nu(x).$$
Further, denoting the operator of the Laguerre derivative by $\beta = DxD$ and its companion $\theta =xDx$ (cf. \cite{Riordan}), where $D$ is the differential operator $D= {d\over dx}$, we calculate them $n$-th power, employing amazing Viskov-type identities \cite{Viskov}

$$\beta^n = \left( DxD\right)^n = D^n x^n D^n,\quad  \theta^n = \left( xDx\right)^n = x^n D^n x^n,\quad   n \in\mathbb{N}_{0}.\eqno(2.6)$$
Equalities  (2.6) can be proved by the method of mathematical  induction. We show how to establish (2.6), using the Mellin transform technique for a class of functions $f$ whose Mellin transforms (1.9) $f^*(s),\ s=\gamma+i\tau$ belong to the Schwartz space as a function of $\tau$. As it is known, this  space  is a topological vector space of functions $\varphi$ such that $ \varphi \in C^\infty (\mathbb{R})$ and $x^m \varphi^{(n)} (x) \to 0, \ |x| \to \infty,\ m,n  \in \mathbb{N}_0$. This means that one can differentiate under the integral sign in (1.11) infinitely many times.  Hence 

$$(\beta^n  f) (x) = \left( DxD\right)^n  f = \frac{1}{2\pi i}  \left( DxD\right)^{n-1}  \int_{\gamma-i\infty}^{\gamma+i\infty}  s^2 f^*(s)  x^{-s-1} ds$$

$$= \frac{1}{2\pi i}  \left( DxD\right)^{n-2}  \int_{\gamma-i\infty}^{\gamma+i\infty}  [ s(s+1)]^2 f^*(s)  x^{-s-2} ds = \dots = \frac{1}{2\pi i}   \int_{\gamma-i\infty}^{\gamma+i\infty}  [ (s)_n]^2 f^*(s)  x^{-s-n} ds, $$
where

$$(s)_n= s(s+1)\dots (s+n-1)= {\Gamma(s+n) \over \Gamma(s)}\eqno(2.7)$$
is the Pochhammer symbol \cite{Bateman}.  On the other hand,

$$ \left( D^n x^n D^n \right)  f = \frac{(-1)^n }{2\pi i}   D^n  \int_{\gamma-i\infty}^{\gamma+i\infty}  (s)_n  f^*(s)  x^{-s} ds
= \frac{1 }{2\pi i}    \int_{\gamma-i\infty}^{\gamma+i\infty}  [ (s)_n ]^2  f^*(s)  x^{-s-n} ds,$$
which proves the first identity in (2.6).  Analogously,

$$ (\theta^n  f) (x) = \left( xDx\right)^n  f = \frac{1}{2\pi i}  \left( xDx\right)^{n-1}  \int_{\gamma-i\infty}^{\gamma+i\infty}  (1-s) f^*(s)  x^{1-s} ds$$

$$= \dots =  \frac{1}{2\pi i}   \int_{\gamma-i\infty}^{\gamma+i\infty}  (1-s)_n f^*(s)  x^{n-s} ds=   \frac{x^n }{2\pi i} D^n   \int_{\gamma-i\infty}^{\gamma+i\infty}  f^*(s)  x^{n-s} ds  = \left( x^n D^n x^n\right)  f. $$
This proves the second identity in (2.6).   In particular, we easily find the values

$$ (\beta^n  \rho_0) (x)= \left( DxD\right)^n  \rho_0 = \rho_0(x), \quad   (\beta^n  \rho_1) (x) = \left( DxD\right)^n  \rho_1 = \rho_1(x) - n \rho_0(x),\quad   n \in\mathbb{N}_{0},\eqno(2.8)$$
 
 $$ (\theta^n  1) (x) = \left( xDx\right)^n  1= n! x^n,\quad  (\theta^n  x^k) (x)=  \left( xDx\right)^n  x^k = {(n+k)!\over k!} x^{n+k},  \quad  n, k  \in\mathbb{N}_{0}.\eqno(2.9)$$
The quotient of the scaled Macdonald functions $\rho_\nu, \rho_{\nu+1}$ is given by the important Ismail integral representation \cite{Ismail}

$${\rho_\nu(x) \over \rho_{\nu+1}(x) } = {1\over \pi^2} \int_0^\infty {y^{-1} dy \over (x+y) \left[ J_{\nu+1}^2 (2\sqrt y)+ Y^2_{\nu+1}(2\sqrt y) \right] },\eqno(2.10)$$
where $J_\nu(z), Y_\nu(z)$ are Bessel functions of the first and second kind, respectively \cite{Bateman}. Another interesting integral representation for the scaled Macdonald function $\rho_\nu$ is given via relation (2.19.4.13) in \cite{PrudnikovMarichev}, Vol. II in terms of the associated Laguerre polynomials. Namely, we have

$${(-1)^n x^n\over n!}\  \rho_\nu(x)=   \int_0^\infty t^{\nu+n -1} e^{-t - x/t}  L_n^\nu(t) dt,\quad    n \in\mathbb{N}_{0}.\eqno(2.11)$$
Meanwhile, important property for  the scaled Macdonald functions can be obtained in terms of the Riemann-Liouville fractional integral \cite{YaL}

$$ \left( I_{-}^\nu  f \right) (x)  = {1\over \Gamma(\nu)} \int_x^\infty (t-x)^{\nu-1} f(t) dt.\eqno(2.12)$$
In fact, appealing to relation (2.16.3.8) in \cite{Bateman}, Vol. II

$$2^{\alpha -1} x^{\alpha+\nu} \Gamma(\alpha) K_{\nu+\alpha}(x) = \int_x^\infty t^{1+\nu} (t^2-x^2)^{\alpha-1} K_\nu(t) dt,\eqno(2.13)$$  
 making simple changes of variables and letting $\alpha =0$, we derive the formula 
$$\rho_\nu(x)= \left( I_{-}^\nu \rho_0 \right) (x).\eqno(2.14)$$
Moreover, the index law for fractional integrals immediately implies

$$ \rho_{\nu+\mu} (x)= \left( I_{-}^\nu \rho_\mu \right) (x)=   \left( I_{-}^\mu \rho_\nu \right) (x).\eqno(2.15)$$
The corresponding definition of the fractional derivative presumes the relation $ D^\mu_{-}= - D  I_{-}^{1-\mu}$.   Hence for the ordinary $n$-th derivative of $\rho_\nu$ we get

$$D^n \rho_\nu(x)= (-1)^n \rho_{\nu-n} (x),\quad n \in \mathbb{N}_0.\eqno(2.16)$$
Another way to get this formula is to differentiate $n$-times  the integral (1.4), to use the definition of the Pochhammer symbol (2.7) and to make a simple change of variables.   

In the meantime, the Mellin-Barnes integral (1.4) and reduction formula for the gamma-function yield 

$$\rho_{\nu+1} (x) = \frac{1}{2\pi i} \int_{\gamma-i\infty}^{\gamma+i\infty} \Gamma(\nu+s+1) \Gamma (s) x^{-s} ds$$

$$=  \frac{1}{2\pi i} \int_{\gamma-i\infty}^{\gamma+i\infty} \Gamma(\nu+s) (\nu+s) \Gamma (s) x^{-s} ds=  \nu \rho_\nu(x)$$

$$+  \frac{1}{2\pi i} \int_{\gamma-i\infty}^{\gamma+i\infty} \Gamma(\nu+s) \Gamma (s+1) x^{-s} ds$$

$$=   \nu \rho_\nu(x)+ x \rho_{\nu-1} (x).$$
Hence we deduce the following recurrence relation for the scaled Macdonald functions

$$\rho_{\nu+1} (x) =    \nu \rho_\nu(x)+ x \rho_{\nu-1} (x),\quad \nu \in \mathbb{R}.\eqno(2.17)$$
In the operator form it can be written as follows

$$\rho_{\nu+1} (x) = \left( \nu - xD \right) \rho_\nu(x),\eqno(2.18)$$
and more generally

$$\rho_{\nu+n} (x) = \prod_{k=0}^{n-1}  \left( \nu+n-k-1 - xD \right) \rho_\nu(x),\quad n \in \mathbb{N}_0.\eqno(2.19)$$
Further, recalling the definition of the operator $\theta$, identities (2.6)  and Rodrigues formula for the associated Laguerre polynomials,  we obtain

$$\theta^n \{ x^\nu e^{-x} \} =  n! x^{n+\nu} e^{-x} L_n^\nu(x),\  n \in \mathbb{N}_0.\eqno(2.20)$$
This formula permits to derive an integral representation for  the product $\rho_\nu f_n$, where $f_n$ is an arbitrary polynomial of degree $n$

$$f_n(x)= \sum_{k=0}^n f_{n,k} x^k.$$
In fact, considering the operator equality and using (2.20), we write

$$f_n(-\theta) \left\{ x^\nu e^{-x}\right\} =  x^\nu e^{-x}  \sum_{k=0}^n f_{n,k} (-1)^k k! x^k L_k^\nu(x) = x^\nu e^{-x} q^\nu_{2n}(x),\eqno(2.21)$$
where 

$$  q^\nu_{2n}(x)= \sum_{k=0}^n f_{n,k} (-1)^k k! x^k L_k^\nu(x)\eqno(2.22)$$
will be called the associated polynomial of degree $2n$.  Then, integrating by parts in the following integral and eliminating the integrated terms, we find

$$\int_0^\infty t^{-1} e^{-x/t}  f_n(-\theta) \left\{ t^\nu e^{-t} \right\} dt  = \int_0^\infty f_n(\theta) \left\{ t^{-1} e^{-x/t} \right\} 
 t^\nu e^{-t} dt .$$
Meanwhile,

$$\theta^k \left\{ t^{-1} e^{-x/t} \right\}  = \left( t D t\right)^k \left\{ t^{-1} e^{-x/t} \right\} = x^k  t^{-1}  e^{-x/t}.$$
Hence, appealing to (2.3) and (2.22), we establish the following integral representation of an arbitrary polynomial $f_n$ in terms of its associated polynomial $q^\nu_{2n}$

$$f_n(x)= {1\over \rho_\nu(x) }  \int_0^\infty t^{\nu-1} e^{-t -x/t  }  q^\nu_{2n}(t) dt.\eqno(2.23)$$ 

The following lemma gives  the so-called linear polynomial independence  of the  scaled Macdonald functions.  Precisely, we have
 
{ \bf Lemma 1}. {\it Let $n, m \in \mathbb{N},\nu \ge 0,\   f_n,\ g_{m}$ be polynomials of  degree at most $n,\ m$, respectively.  Let 
 
 $$f_n(x) \rho_\nu(x)+ g_{m} (x) \rho_{\nu+1} (x) = 0\eqno(2.24)$$
 for all $x >0$. Then $f_n \equiv 0, \ g_{m} \equiv 0.$ }
 
 \begin{proof} The proof will be based on the Ismail integral representation (2.10) of the quotient $\rho_\nu/ \rho_{\nu+1}$.  In fact,  since $\rho_{\nu+1}  > 0$, we divide (2.24)  by  $\rho_{\nu+1} $ and  then differentiate $m+1$ times the obtained equality. Thus we arrive at  the relation
 
 $${d^{m+1}\over dx^{m+1} } \left[ f_n(x) {\rho_\nu(x)\over \rho_{\nu+1} (x) }\right]  = 0,\quad  x >0.\eqno(2.25)$$
 Meanwhile, the integral representation (2.10) says

$$ {\rho_{\nu}(x) \over \rho_{\nu+1}(x) } = {1\over \pi^2} \int_0^\infty { s^{-1} ds \over (x+s)( J_{\nu+1}^2(2\sqrt s) + 
Y_{\nu+1}^2(2\sqrt s) )}$$

$$ = {1\over \pi^2} \int_0^\infty e^{-xy } dy \int_0^\infty { e^{-sy} \ s^{-1} ds \over  J_{\nu+1}^2(2\sqrt s) + 
Y_{\nu+1}^2(2\sqrt s)},\eqno(2.26)$$
where the interchange of the order of integration is allowed by Fubini theorem, taking into account the asymptotic behavior of Bessel functions at infinity and near zero \cite{Bateman}. Further, assuming that

$$f_n(x)= \sum_{k=0}^n f_{n,k}\  x^k,$$
we substitute it in the left-hand side of (2.25) together with the right-hand side of the latter equality in (2.26). Then,  differentiating  under the integral sign, which is possible via the absolute and uniform convergence, we deduce 

$$  {d^{m+1}\over dx^{m+1} } \left[ f_n(x) {\rho_\nu(x)\over \rho_{\nu+1} (x) }\right]  =  {1\over \pi^2}  {d^{m+1}\over dx^{m+1} }   \sum_{k=0}^n f_{n,k} x^k \int_0^\infty e^{-xy } dy \int_0^\infty { e^{-sy} \ s^{-1} ds \over  J_{\nu+1}^2(2\sqrt s) + 
Y_{\nu+1}^2(2\sqrt s)}$$

$$=  {1\over \pi^2}    \sum_{k=0}^n f_{n,k} (-1)^k  {d^{m+1}\over dx^{m+1} }  \int_0^\infty  {d^k\over dy^k} \left[ e^{-xy } \right] dy \int_0^\infty { e^{-sy} \ s^{-1} ds \over  J_{\nu+1}^2(2\sqrt s) + Y_{\nu+1}^2(2\sqrt s)}$$

$$=  {1\over \pi^2}    \sum_{k=0}^n f_{n,k} (-1)^k   \int_0^\infty  {\partial^{k+m+1} \over \partial y^k \partial x^{m+1} } \left[ e^{-xy } \right] dy \int_0^\infty { e^{-sy} \ s^{-1} ds \over  J_{\nu+1}^2(2\sqrt s) + Y_{\nu+1}^2(2\sqrt s)}$$

 $$=  {1\over \pi^2}    \sum_{k=0}^n f_{n,k} (-1)^{k+m+1}    \int_0^\infty  {d^k\over dy^k} \left[  y^{m+1}  e^{-xy } \right] dy \int_0^\infty { e^{-sy} \ s^{-1} ds \over  J_{\nu+1}^2(2\sqrt s) + Y_{\nu+1}^2(2\sqrt s)}.$$
Now, integrating $k$ times by parts in the outer integral with respect to $y$ on the right-hand side of the latter equality,  and then differentiating under the  integral sign in the inner integral with respect to $s$ owing to the same arguments, we get, combining with (2.25)

$$ {1\over \pi^2}    \sum_{k=0}^n f_{n,k} (-1)^{k+m+1}    \int_0^\infty  {d^k\over dy^k} \left[  y^{m+1}  e^{-xy } \right] dy \int_0^\infty { e^{-sy} \ s^{-1} ds \over  J_{\nu+1}^2(2\sqrt s) + Y_{\nu+1}^2(2\sqrt s)}$$

 $$ = {1\over \pi^2}     \int_0^\infty   y^{m+1}  e^{-xy }  \int_0^\infty { e^{-sy} \ s^{-1}  \over  J_{\nu+1}^2(2\sqrt s) + Y_{\nu+1}^2(2\sqrt s)} 
 \left(    \sum_{k=0}^n f_{n,k} (-1)^{k+m+1}  s^k \right) ds = 0,\ x >0.\eqno(2.27)$$
Consequently, cancelling twice  the Laplace transform (1.6) via its injectivity for integrable continuous functions \cite{Tit},  and taking into account the positivity of the function 

$$ { s^{-1}  \over  J_{\nu+1}^2(2\sqrt s) + Y_{\nu+1}^2(2\sqrt s)} $$
on $\mathbb{R}_+$,  we conclude that

$$ \sum_{k=0}^n f_{n,k} (-1)^{k}  s^k  \equiv 0,\quad s >0.$$
Hence $f_{n,k} =0,\ k=0,\dots, n$  and therefore $f_n \equiv 0.$  Returning to the original equality (2.24), we find immediately that $g_{m} \equiv 0.$ Lemma 1  is proved. 
 
  \end{proof}

 Let $\alpha \in \mathbb{R}$ and 
 
 $$  S^{\nu,\alpha}_n(x)= {d^n\over dx^n} \left[ x^{n+\alpha} \rho_\nu(x)\right],\quad n \in \mathbb{N}_0.\eqno(2.28)$$
According to \cite{AsscheYakubov2000},  the sequence of functions   $\left(S^{\nu,\alpha}_n\right)_{n\in\mathbb{N}_0}$ generates multiple orthogonal polynomials related to the scaled Macdonald functions $\rho_\nu,\ \rho_{\nu+1}$. In order to obtain an integral representation for functions  $  S^{\nu,\alpha}_n$, we employ again (1.4), Parseval equality (1.10) for the Mellin transform and the Mellin-Barnes integral representation for the associated Laguerre polynomials (see relation  (8.4.33.3) in \cite{PrudnikovMarichev}, Vol. III ).  Then, motivating the differentiation under the integral sign by the absolute and uniform  convergence and using the reflection formula for the gamma-function, we obtain the following chain of equalities

$${d^n\over dx^n} \left[ x^{n+\alpha} \rho_\nu(x)\right] = {1\over 2\pi i} {d^n\over dx^n} \int_{\gamma-i\infty}^{\gamma+i\infty} \Gamma (s+n+\alpha )\Gamma(s+n+ \nu+\alpha) x^{-s} ds$$

$$= {(-1)^n \over 2\pi i} \int_{\gamma-i\infty}^{\gamma+i\infty} \Gamma (s+n+\alpha )\Gamma(s+n+ \nu+\alpha) (s)_n \  x^{-s-n} ds$$

$$ = {(-1)^n \over 2\pi i} \int_{\gamma+n-i\infty}^{\gamma+n +i\infty} \Gamma (s+\alpha )\Gamma(s+ \nu+\alpha) {\Gamma(s)\over \Gamma(s-n)}  \  x^{-s} ds$$

$$= {(-1)^n x^\alpha \over 2\pi i} \int_{\gamma+n+\alpha -i\infty}^{\gamma+n+\alpha  +i\infty} \Gamma (s )\Gamma(s+ \nu) {\Gamma(s-\alpha)\over \Gamma(s-\alpha -n)}  \  x^{-s} ds$$

$$= { x^\alpha \over 2\pi i} \int_{\gamma -i\infty}^{\gamma  +i\infty} \Gamma (s )\Gamma(s+ \nu) {\Gamma(1+\alpha+n-s)\over \Gamma(1+\alpha-s)}  \  x^{-s} ds$$

$$= { x^{\alpha+\nu}  \over 2\pi i} \int_{\gamma+\nu  -i\infty}^{\gamma +\nu  +i\infty} \Gamma (s-\nu )\Gamma(s) {\Gamma(1+\alpha+\nu+ n-s)\over \Gamma(1+\alpha+\nu -s)}  \  x^{-s} ds$$

$$=   x^{\alpha+\nu} n! \int_0^\infty  e^{-t-x/t} \left({x\over t} \right)^{-\nu}  L_n^{\nu+\alpha} (t) {dt\over t} .$$
Thus, combining with (2.28), we established the following integral representation for    $S^{\nu,\alpha}_n(x)$ 

$$S^{\nu,\alpha}_n(x) =   x^{\alpha} n! \int_0^\infty  e^{-t-x/t} t^{\nu-1}  L_n^{\nu+\alpha} (t) dt,\quad  x >0. \eqno(2.29)$$ 
Now, employing recurrence relations and differential properties for the associated Laguerre polynomials \cite{Chi}, making integration by parts in (2.29) and differentiating with respect to $x$ under the integral sign by virtue of the absolute and uniform convergence with respect to $ x \ge x_0 >0$, we will deduce the corresponding relations for the sequence $S^{\nu,\alpha}_n$.  Indeed, we have, for instance, for $\nu >0, \alpha \in \mathbb{R}$

$$S^{\nu+1,\alpha-1}_n(x) =   x^{\alpha-1} n! \int_0^\infty  e^{-t-x/t} t^{\nu}  L_n^{\nu+\alpha} (t) dt$$ 

$$=   x^{\alpha-1} n! \left[ \nu \int_0^\infty  e^{-t-x/t} t^{\nu-1}  L_n^{\nu+\alpha} (t) dt + x  \int_0^\infty  e^{-t-x/t} t^{\nu-2}  L_n^{\nu+\alpha} (t) dt \right.$$ 

$$\left.  -  \int_0^\infty  e^{-t-x/t} t^{\nu}  L_{n-1}^{\nu+\alpha+1} (t) dt \right] =  {\nu \over x} \ S^{\nu,\alpha}_n(x) + {1 \over x} \ S^{\nu-1,\alpha+1}_n(x)
-  {n \over x} \ S^{\nu+1,\alpha}_{n-1}(x).$$
Hence we obtain the identity

$$x S^{\nu+1,\alpha-1}_n(x) = \nu \ S^{\nu,\alpha}_n(x) +  S^{\nu-1,\alpha+1}_n(x) -  n\  S^{\nu+1,\alpha}_{n-1}(x),\ x > 0,\ n \in \mathbb{N}_0.\eqno(2.30)$$
Differentiating (2.29) by $x$, we get

$${d\over dx} \  S^{\nu,\alpha}_n(x)   =   \alpha x^{\alpha-1} n! \int_0^\infty  e^{-t-x/t} t^{\nu-1}  L_n^{\nu+\alpha} (t) dt -  x^{\alpha} n! \int_0^\infty  e^{-t-x/t} t^{\nu-2}  L_n^{\nu+\alpha} (t) dt,$$
or,

$$  x {d\over dx} \  S^{\nu,\alpha}_n(x)   = \alpha  S^{\nu,\alpha}_n(x) -  S^{\nu-1,\alpha+1}_n(x),  \ x > 0,\ n \in \mathbb{N}_0.\eqno(2.31)$$
On the other hand, integrating again   by parts in (2.29) under the same conditions, we find  

$$  S^{\nu,\alpha}_n(x) = {x^\alpha n! \over \nu} \int_0^\infty e^{-t-x/t} t^{\nu} L_n^{\nu+\alpha} (t) dt + {x^\alpha n!  \over \nu} \int_0^\infty e^{-t-x/t} t^{\nu} L_{n-1}^{\nu+\alpha+1} (t) dt $$

$$- {x^{\alpha +1} n! \over \nu} \int_0^\infty e^{-t-x/t} t^{\nu-2 } L_n^{\nu+\alpha} (t) dt  = {1\over \nu} S^{\nu+1,\alpha-1}_n(x) $$

$$+ {1\over \nu} S^{\nu+1,\alpha}_{n-1} (x) - {1\over \nu} S^{\nu-1,\alpha+1}_n(x),$$ 
or,

$$   \nu S^{\nu,\alpha}_n(x) = S^{\nu+1,\alpha-1}_n(x) + S^{\nu+1,\alpha}_{n-1} (x) - S^{\nu-1,\alpha+1}_n(x), \ x > 0,\ n \in \mathbb{N}_0.\eqno(2.32)$$
Combining with (2.30) it gives the following identity

$$(x-1)  S^{\nu+1,\alpha-1}_n(x) =  (1-n) S^{\nu+1,\alpha}_{n-1} (x),\ x > 0,\ n \in \mathbb{N}_0.\eqno(2.33)$$
Meanwhile, from (2.28) and (2.17)  it has 
$$ S^{\nu-1,\alpha+1}_n(x) =  {d^n\over dx^n} \left[ x^{n+\alpha+1} \rho_{\nu-1} (x)\right] $$

$$=  {d^n\over dx^n} \left[ x^{n+\alpha} \left[ \rho_{\nu+1}(x)- \nu \rho_\nu(x) \right] \right] $$

$$= S^{\nu+1,\alpha}_n(x)-  \nu  S^{\nu,\alpha}_n(x).$$
Therefore from (2.32) we have 

$$  S^{\nu+1,\alpha}_n(x) = S^{\nu+1,\alpha}_{n-1} (x) +   S^{\nu+1,\alpha-1}_n(x),\eqno(2.34)$$
and from (2.33) we find

$$(x-1)  S^{\nu+1,\alpha}_n(x)   =  (x-n) S^{\nu+1,\alpha}_{n-1} (x),\ x > 0,\ n \in \mathbb{N}_0.\eqno(2.35)$$
Moreover, recalling again (2.17), we deduce
$${d\over dx}  S^{\nu+1,\alpha}_{n-1} (x)  =  {d^{n}\over dx^{n} } \left[ x^{n+\alpha-1} \rho_{\nu+1} (x)\right]  =  \nu  S^{\nu,\alpha-1}_{n} (x) +   S^{\nu-1,\alpha}_n(x).\eqno(2.36)$$
Finally, employing the 3-term recurrence relation for the associated Laguerre polynomials

$$(n+1)L_{n+1}^{\nu+\alpha} (x)= (2n+1+\nu+\alpha -x)  L_{n}^{\nu+\alpha} (x) - (n+\nu+\alpha) L_{n-1}^{\nu+\alpha} (x),\eqno(2.37)$$
we return to (2.29) to obtain the following identity 

$$S^{\nu,\alpha}_{n+1} (x) =   (2n+1+\nu+\alpha ) S^{\nu,\alpha}_{n} (x)   -  n (n+\nu+\alpha) S^{\nu,\alpha}_{n-1} (x) -  x S^{\nu+1,\alpha-1}_{n} (x), \ x > 0,\ n \in \mathbb{N}_0.\eqno(2.38)$$

\section{ Prudnikov's orthogonal polynomials} 

Our goal in this section is to find an explicit expression for Prudnikov's orthogonal polynomial sequence $\left(P_n^\nu\right)_{ n\in \mathbb{N}_0},\ \nu \ge 0.$ We will do even more, defining the Prudnikov orthogonality (1.5) in a more general setting for the sequence  $\left(P_n^{\nu,\alpha}\right)_{ n\in \mathbb{N}_0},\ \alpha > -1,$

$$\int_0^\infty  P_n^{\nu,\alpha} (x)P_m^{\nu,\alpha} (x)  x^\alpha \rho_\nu(x) dx=  \delta_{n,m} , \quad n,m\in\mathbb{N}_{0}.\eqno(3.1)$$
Here  $P_n^\nu \equiv  P_n^{\nu,0}.$  Writing it in terms of coefficients

$$P_n^{\nu,\alpha} (x)= \sum_{k=0}^n a_{n,k} x^k,\eqno(3.2)$$
we know  that it is of degree exactly $n$ because this sequence is regular, i.e. its leading coefficient $a_{n,n}\equiv a_n \neq 0$ (cf. \cite{Filipa}).  Furthermore, as it follows from the general theory of orthogonal polynomials \cite{Chi}, up to a normalization factor the orthogonality (3.1) is equivalent to the following $n$ conditions

$$\int_0^\infty  P_n^{\nu,\alpha} (x) x^{m+\alpha}  \rho_\nu(x) dx = 0,\quad m= 0,1, \dots, n-1.\eqno(3.3)$$
Moreover, the sequence $ \left(P_n^{\nu,\alpha}\right)_{ n\in \mathbb{N}_0}$ satisfies the 3-term recurrence relation in the form

$$x  P_n^{\nu,\alpha} (x) = A_{n+1}  P_{n+1}^{\nu,\alpha} (x) + B_n  P_n^{\nu,\alpha} (x) + A_n  P_{n-1}^{\nu,\alpha} (x),\eqno(3.4)$$
where $P_{-1}^{\nu,\alpha} (x)\equiv 0$ and

$$A_{n+1}= {a_n\over a_{n+1} }, \quad B_{n}= {b_n\over a_n} -  {b_{n+1}\over a_{n+1}},\quad  b_n\equiv a_{n,n-1}.\eqno(3.5)$$
The associated polynomial sequence to $\left(P_n^{\nu,\alpha}\right)_{ n\in \mathbb{N}_0}$ (cf. (2.22)), which will be used in the sequel, has the form

$$Q_{2n}(x)= \sum_{k=0}^n a_{n,k} (-1)^k k! x^k L_k^\nu(x).\eqno(3.6)$$
As it follows from the orthogonality (3.1)

$$\int_0^\infty  \left[ P_n^{\nu,\alpha} (x) \right]^2 x^\alpha \rho_\nu(x) dx =1.$$
However, using properties of the scaled Macdonald functions from the previous section one can calculate the following values

$$\int_0^\infty  \left[ P_n^{\nu,\alpha} (x) \right]^2  x^\alpha \rho_{\nu+1} (x) dx.$$
In fact, appealing to (3.1), (3.2), (3.6), (2.17)  and integrating by parts, we derive

$$\int_0^\infty  \left[ P_n^{\nu,\alpha} (x) \right]^2 x^{\alpha} \rho_{\nu+1} (x) dx= \nu+  \int_0^\infty  \left[ P_n^{\nu,\alpha} (x) \right]^2  x^{\alpha+1}  \rho_{\nu-1} (x) dx$$

$$= \nu +\alpha+ 1+ 2 \int_0^\infty   P_n^{\nu,\alpha} (x) {d\over dx} \left[ P_n^{\nu,\alpha} (x) \right]  x^{\alpha+1} \rho_{\nu} (x) dx = 2n+1+\nu+\alpha$$
since

$$ \int_0^\infty   P_n^{\nu+\alpha} (x)  x^{n+\alpha}  \rho_\nu(x) dx = {1\over a_n}.\eqno(3.7)$$
Therefore we find the formula

$$\int_0^\infty  \left[ P_n^{\nu+\alpha} (x) \right]^2  x^\alpha \rho_{\nu+1} (x) dx = 2n+1+\nu+\alpha.\eqno(3.8)$$
In the meantime, taking the corresponding integral representation (2.11) for the product $x^m \rho_\nu(x)$, we substitute its right-hand side in (3.3) and change the order of integration by Fubini's theorem. Thus we obtain

$$ \int_0^\infty t^{\nu+m -1} e^{-t}  L_m^\nu(t) \int_0^\infty  P_n^{\nu,\alpha} (x) e^{-x/t}  x^\alpha dx dt=0, \quad m= 0,1, \dots, n-1.\eqno(3.9)$$
But the inner integral with respect to $x$ can be treated, involving the differential operator $\theta$ (see (2.6)).  Indeed, using (3.2), we have

$$ {1\over t} \int_0^\infty  P_n^{\nu,\alpha} (x) e^{-x/t}  x^\alpha dx =  \sum_{k=0}^n a_{n,k}  \theta^k \left\{ {1\over t } \int_0^\infty  e^{-x/t} x^\alpha dx \right\}$$

$$ =  \Gamma(1+\alpha) \sum_{k=0}^n a_{n,k}  \theta^k \left\{ t^\alpha \right\} =  \Gamma(1+\alpha) P_n^{\nu,\alpha}  (\theta) \{t^\alpha \}.$$
Moreover, the Rodrigues formula for the associated Laguerre polynomials and Viskov type identity (2.6) for the operator $\theta$ imply

$$ t^{\nu+m} e^{-t}  L_m^\nu(t) = {1\over n!} \theta^m \left\{ t^\nu e^{-t} \right\}.$$
Substituting these values in (3.9), it becomes

$$ \int_0^\infty   \theta^m \left\{ t^\nu e^{-t} \right\} P_n^{\nu,\alpha}  (\theta) \{t^\alpha \} dt=0, \quad m= 0,1, \dots, n-1.$$
After $m$ times  integration by parts in the latter integral, we end up with the following orthogonality conditions 

$$ \int_0^\infty  t^\nu e^{-t}   \theta^m P_n^{\nu,\alpha}  (\theta) \{t^\alpha \} dt=0, \quad m= 0,1, \dots, n-1.\eqno(3.10)$$
Analogously, the orthogonality (3.1) is equivalent to the equality

$$ \int_0^\infty  t^\nu e^{-t}   P^{\nu,\alpha}_m (\theta) P_n^{\nu,\alpha}  (\theta) \{t^\alpha \} dt= {\delta_{m,n}\over \Gamma(1+\alpha)},\quad  \alpha > -1.\eqno(3.11)$$

{\bf Definition 2}. {\it The orthogonality $(3.11)$ is called the composition orthogonality of the sequence $\left(P_n^{\nu,\alpha}\right)_{ n\in \mathbb{N}_0} $ in the sense of Laguerre.}

Thus we proved the following theorem.

{\bf Theorem 1.} {\it The Prudnikov orthogonality $(3.1)$ is equivalent to the composition orthogonality $(3.11)$ in the  sense of Laguerre, i.e. Prudnikov's orthogonal polynomials are the associated Laguerre polynomials in the sense of composition orthogonality $(3.11)$.}

Meanwhile, in terms of the associated polynomial (3.6) the orthogonality conditions (3.10) can be rewritten, using the commutativity property

$$  \theta^m P_n^{\nu,\alpha} (\theta) \{t^\alpha\}  =   P_n^{\nu,\alpha}  (\theta) \theta^m \{t^\alpha\} $$
and the Rodrigues formula for the associated Laguerre polynomials. Then, integrating by parts an appropriate number of times and  taking into account (2.9), we get

$$ 0=  \int_0^\infty  t^\nu e^{-t}   \theta^m P_n^{\nu,\alpha}  (\theta) \{t^\alpha\} dt =  \int_0^\infty  t^\nu e^{-t}   P_n^{\nu,\alpha}  (\theta) \theta^m \{t^\alpha\} dt$$

$$=  (1+\alpha)_m  \int_0^\infty P^{\nu,\alpha}_n\left( -\theta\right) \left\{  t^\nu e^{-t} \right\} t^{m+\alpha}   dt =  (1+\alpha)_m  \int_0^\infty  t^{\nu+ \alpha+ m}  e^{-t} Q_{2n} (t) dt,$$
or, finally,

$$ \int_0^\infty  t^{\nu+\alpha+m}  e^{-t}  Q_{2n} (t) dt = 0,  \quad m= 0,1, \dots, n-1.\eqno(3.12)$$
On the other hand, developing the polynomial $ Q_{2n} (t)$ in terms of the associated Laguerre polynomials $L_n^{\nu+\alpha}(x)$, we find

$$ Q_{2n} (x) = \sum_{j=0}^{2n} c_{n,j} L_j^{\nu+\alpha}(x),\eqno(3.13)$$
where

$$c_{n,k}=  {k!\over \Gamma(k+\nu+\alpha+1)}  \int_0^\infty  t^{\nu+\alpha}  e^{-t}  Q_{2n} (t) L_k^{\nu+\alpha} (t) dt\eqno(3.14)$$ 
and orthogonality conditions (3.12) immediately imply that

$$c_{n,j} = 0, \quad j = 0,1, \dots, n-1.\eqno(3.15)$$
Therefore, the expansion (3.13) becomes

$$ Q_{2n} (x) = \sum_{j=n}^{2n} c_{n,j} L_j^{\nu+\alpha} (x).\eqno(3.16)$$
In the meantime,  expanding  $(-1)^m  m! x^m\  L_m^\nu(x)$ in terms of the associated Laguerre polynomials $L_k^{\nu+\alpha} (x)$ as well, we obtain

$$ (-1)^m  m! x^m\  L_m^\nu(x) = \sum_{k=0}^{2m} d_{m,k} L_k^{\nu+\alpha} (x),\eqno(3.17)$$
where  coefficients $d_{m,k}$ are calculated accordingly by the formula (see relation (2.19.14.8)  in \cite{PrudnikovMarichev}, Vol. II)

$$d_{m,k}= {(-1)^m\  m! \ k! \over \Gamma (k+\nu+\alpha+1)} \int_0^\infty  e^{-t}  t^{\nu+\alpha +m} L_m^\nu(t) \   L_k^{\nu+\alpha} (t) dt$$

$$=  {(-1)^{m+k}\  m! \over (m-k)!}  \  (1+\nu)_m\ (\nu+\alpha+1+k)_{m-k}   \ 
{}_3F_2 \left( -m,\ \nu+\alpha +m+1,\ m+1;\ 1+\nu,\ m+1-k;\ 1\right), \eqno(3.18)$$
where ${}_3F_2(a,b,c; d,e;z)$ is the generalized hypergeometric function \cite{PrudnikovMarichev}, Vol. III.   It is easily seen from the orthogonality of the associated Laguerre polynomials $L_k^{\nu+\alpha} (x)$ that 

$$d_{m,k} =0,\quad  k > 2m.\eqno(3.19)$$
Moreover,  the associated polynomial (3.6) $Q_{2n}$ has the representation

$$Q_{2n}(x) =   \sum_{m=0}^n  \ a_{n,m}   \sum_{k=0}^{2m} d_{m,k} L_k^{\nu+\alpha} (x) $$
$$=   \sum_{m=0}^n  \ a_{n,m}  \left[  \sum_{k=0}^{m} d_{m, 2k} L_{2k}^{\nu+\alpha} (x) +  \sum_{k=0}^{m-1} d_{m, 2k+1} L_{2k+1}^{\nu+\alpha} (x) \right] $$

$$=    \sum_{k=0}^{n}  L_{2k}^{\nu+\alpha} (x)  \left( \sum_{m= k}^n  \ a_{n,m}  \ d_{m, 2k} \right)  +  \sum_{k=0}^{n-1}  L_{2k+1}^{\nu+\alpha} (x)  \left( \sum_{m= k}^{n-1}   \ a_{n,m+1}  \ d_{m+1, 2k+1} \right) .$$

{\bf Lemma 2}. {\it Coefficients $d_{m,k},\ m,k \in \mathbb{N}_0,$ satisfy the following recurrence relation}

$$ d_{m+1,k}= - m k(k-1) (m+\nu) d_{m-1,k-2} +  m k (m+\nu) (1+2\alpha+3k+2(m+\nu)) d_{m-1,k-1} $$

$$- m(m+\nu)  \left(\alpha^2+ 3k^2+ (\nu+1)(2(1+m)+\nu)+ \alpha(3+4k+ 2(m+\nu) ) + k (5+ 4(m+\nu)) \right) d_{m-1,k} $$

$$+ m(m+\nu)(1+\alpha+k+\nu)(2(1+m) +\alpha+k +\nu) d_{m-1,k +1} + k(k-1) d_{m,k-2}$$ 

$$- k(3k+2\alpha+\nu)  d_{m,k-1} +  \left( (1+\alpha)(1+\alpha+4k)+ 3(k^2-m^2) + 2m(k-1)\right.$$

$$\left. \left. + \nu(1+\alpha+3k-m) \right) d_{m,k} - (1+\alpha+k+\nu)(2(1+m) + \alpha+k+\nu) \right) d_{m,k+1}.\eqno(3.20)$$

\begin{proof}  In fact, recalling the 3-term recurrence relation (2.37) for the associated Laguerre polynomials and,  as its direct consequence,  the following equality

$$x L_{n}^{\nu+\alpha+1} (x) = (n+\nu+\alpha) L_{n-1}^{\nu+\alpha} (x) - (n-x) L_{n}^{\nu+\alpha} (x),\eqno(3.21)$$
we derive  from  (3.18) via integration  by parts  

$$ d_{m+1,k}= {(-1)^{m+1}\  (m+1)! \ k! \over \Gamma (k+\nu+\alpha+1)} \int_0^\infty  e^{-t}  t^{\nu+\alpha +m+1} L_{m+1}^\nu(t) \   L_k^{\nu+\alpha} (t) dt$$

$$ = {(-1)^{m+1}\  (m+1)! \ k!  (\nu+\alpha +m+1 )\over \Gamma (k+\nu+\alpha+1)} \int_0^\infty  e^{-t}  t^{\nu+\alpha +m} L_{m+1}^\nu(t) \   L_k^{\nu+\alpha} (t) dt$$

$$ + {(-1)^{m}\  (m+1)! \ k!  \over \Gamma (k+\nu+\alpha+1)} \int_0^\infty  e^{-t}  t^{\nu+\alpha +m+1} L_{m}^{\nu+1}(t) \   L_k^{\nu+\alpha} (t) dt$$

$$ + {(-1)^{m}\  (m+1)! \ k!  \over \Gamma (k+\nu+\alpha+1)} \int_0^\infty  e^{-t}  t^{\nu+\alpha +m+1} L_{m+1}^{\nu}(t) \   L_{k-1}^{\nu+\alpha+1} (t) dt$$

$$ = {(-1)^{m+1}\  m! \ k!  (\nu+\alpha +m+1 )\over \Gamma (k+\nu+\alpha+1)} \int_0^\infty  e^{-t}  t^{\nu+\alpha +m} \left[ (2m+1+\nu -t)  L_{m}^{\nu} (t) - (m+\nu) L_{m-1}^{\nu} (t) \right] \   L_k^{\nu+\alpha} (t) dt$$

$$ + {(-1)^{m}\  (m+1)! \ k!  \over \Gamma (k+\nu+\alpha+1)} \int_0^\infty  e^{-t}  t^{\nu+\alpha +m} \left[ (m+\nu) L_{m-1}^{\nu} (t) - (m-t) L_{m}^{\nu} (t) \right]   L_k^{\nu+\alpha} (t) dt$$

$$ + {(-1)^{m}\  m! \ k!  \over \Gamma (k+\nu+\alpha+1)} \int_0^\infty  e^{-t}  t^{\nu+\alpha +m} \left[ (2m+1+\nu -t)  L_{m}^{\nu} (t) - (m+\nu) L_{m-1}^{\nu} (t) \right] $$

$$\times  \left[ (k+\nu+\alpha-1) L_{k-2}^{\nu+\alpha} (t) - (k-t-1) L_{k-1}^{\nu+\alpha} (t) \right]  dt$$

$$ =  - ( \nu+\alpha +m+1 ) (2m+1+\nu) d_{m,k} +  {(-1)^{m}\  m! \ k!  (\nu+\alpha +m+1 )\over \Gamma (k+\nu+\alpha+1)} \int_0^\infty  e^{-t}  t^{\nu+\alpha +m} L_{m}^{\nu} (t) \left[ (2k+1+\nu+\alpha )  L_{k}^{\nu+\alpha} (t)\right.$$

$$\left.  - (k+\nu+\alpha) L_{k-1}^{\nu+\alpha} (t) -(k+1) L_{k+1}^{\nu+\alpha} (t) \right] dt$$

$$+  {(-1)^{m}\  m! \ k!  (\nu+\alpha +m+1 ) (m+\nu) \over \Gamma (k+\nu+\alpha+1)} \int_0^\infty  e^{-t}  t^{\nu+\alpha +m-1} L_{m-1}^{\nu} (t) \left[ (2k+1+\nu+\alpha )  L_{k}^{\nu+\alpha} (t)\right.$$

$$\left.  - (k+\nu+\alpha) L_{k-1}^{\nu+\alpha} (t) -(k+1) L_{k+1}^{\nu+\alpha} (t) \right] dt$$

$$- m(m+1) d_{m,k} +  {(-1)^{m}\  (m+1)! \ k!  (m+\nu) \over \Gamma (k+\nu+\alpha+1)} \int_0^\infty  e^{-t}  t^{\nu+\alpha +m-1}  L_{m-1}^{\nu} (t) \left[  (2k+1+\nu+\alpha )  L_{k}^{\nu+\alpha} (t)\right.$$

$$\left.  - (k+\nu+\alpha) L_{k-1}^{\nu+\alpha} (t) -(k+1) L_{k+1}^{\nu+\alpha} (t) \right] dt + {(-1)^{m}\  (m+1)! \ k!  \over \Gamma (k+\nu+\alpha+1)} \int_0^\infty  e^{-t}  t^{\nu+\alpha +m}  L_{m}^{\nu} (t) $$

$$ \times \left[ (2k+1+\nu+\alpha )  L_{k}^{\nu+\alpha} (t)  - (k+\nu+\alpha) L_{k-1}^{\nu+\alpha} (t) -(k+1) L_{k+1}^{\nu+\alpha} (t) \right] dt$$

$$+   {k(k-1) (2m+1+\nu) \over  k+\nu+\alpha}  \left[ d_{m,k-2}  -   d_{m,k-1} \right] $$

$$-  {(-1)^{m}\  m! \ k!  ( k+\nu+\alpha-1)  \over \Gamma (k+\nu+\alpha+1)} \int_0^\infty  e^{-t}  t^{\nu+\alpha +m}  L_{m}^{\nu} (t) \left[  (2k-3+\nu+\alpha )  L_{k-2}^{\nu+\alpha} (t) -  (k+\nu+\alpha-2) L_{k-3}^{\nu+\alpha} (t) \right.$$

$$\left. - (k-1)  L_{k-1}^{\nu+\alpha} (t) \right] dt$$

$$+  {(-1)^{m}\  m! \ k!  ( 2m+\nu+k)  \over \Gamma (k+\nu+\alpha+1)} \int_0^\infty  e^{-t}  t^{\nu+\alpha +m}  L_{m}^{\nu} (t) \left[  (2k-1+\nu+\alpha )  L_{k-1}^{\nu+\alpha} (t) -  (k+\nu+\alpha-1) L_{k-2}^{\nu+\alpha} (t) \right.$$

$$\left. - k  L_{k}^{\nu+\alpha} (t) \right] dt$$

$$-  {(-1)^{m}\  m! \ k!  ( k+\nu+\alpha-1) (m+\nu)  \over \Gamma (k+\nu+\alpha+1)} \int_0^\infty  e^{-t}  t^{\nu+\alpha +m-1}  L_{m-1}^{\nu} (t) \left[  (2k-3+\nu+\alpha )  L_{k-2}^{\nu+\alpha} (t) -  (k+\nu+\alpha-2) L_{k-3}^{\nu+\alpha} (t) \right.$$

$$\left. - (k-1)  L_{k-1}^{\nu+\alpha} (t) \right] dt$$

$$+  {(-1)^{m}\  m! \ k!  ( m+\nu) (k-1) \over \Gamma (k+\nu+\alpha+1)} \int_0^\infty  e^{-t}  t^{\nu+\alpha +m-1}  L_{m-1}^{\nu} (t) \left[  (2k-1+\nu+\alpha )  L_{k-1}^{\nu+\alpha} (t) -  (k+\nu+\alpha-1) L_{k-2}^{\nu+\alpha} (t) \right.$$

$$\left. - k  L_{k}^{\nu+\alpha} (t) \right] dt$$

$$-  {(-1)^{m}\  m! \ k! (m+\nu)   \over \Gamma (k+\nu+\alpha+1)} \int_0^\infty  e^{-t}  t^{\nu+\alpha +m-1}  L_{m-1}^{\nu} (t) \left[  (2k-1+\nu+\alpha ) \left[  (2k-1+\nu+\alpha )  L_{k-1}^{\nu+\alpha} (t) -  (k+\nu+\alpha-1) L_{k-2}^{\nu+\alpha} (t) \right.\right.$$

$$\left.\left.  - k  L_{k}^{\nu+\alpha} (t) \right]   -  (k+\nu+\alpha-1) \left[  (2k-3+\nu+\alpha )  L_{k-2}^{\nu+\alpha} (t) -  (k+\nu+\alpha-2) L_{k-3}^{\nu+\alpha} (t) \right.\right.$$

$$\left.\left.  - (k-1)  L_{k-1}^{\nu+\alpha} (t) \right]  - k   \left[  (2k+1+\nu+\alpha )  L_{k}^{\nu+\alpha} (t) -  (k+\nu+\alpha) L_{k-1}^{\nu+\alpha} (t) \right.\right.$$

$$\left.\left.  - (k+1)  L_{k+1}^{\nu+\alpha} (t) \right]  \right] dt$$

$$-  {(-1)^{m}\  m! \ k!   \over \Gamma (k+\nu+\alpha+1)} \int_0^\infty  e^{-t}  t^{\nu+\alpha +m}  L_{m}^{\nu} (t) \left[  (2k-1+\nu+\alpha ) \left[  (2k-1+\nu+\alpha )  L_{k-1}^{\nu+\alpha} (t) -  (k+\nu+\alpha-1) L_{k-2}^{\nu+\alpha} (t) \right.\right.$$

$$\left.\left.  - k  L_{k}^{\nu+\alpha} (t) \right]   -  (k+\nu+\alpha-1) \left[  (2k-3+\nu+\alpha )  L_{k-2}^{\nu+\alpha} (t) -  (k+\nu+\alpha-2) L_{k-3}^{\nu+\alpha} (t) \right.\right.$$

$$\left.\left.  - (k-1)  L_{k-1}^{\nu+\alpha} (t) \right]  - k   \left[  (2k+1+\nu+\alpha )  L_{k}^{\nu+\alpha} (t) -  (k+\nu+\alpha) L_{k-1}^{\nu+\alpha} (t) \right.\right.$$

$$\left.\left.  - (k+1)  L_{k+1}^{\nu+\alpha} (t) \right]  \right] dt$$

$$=     ( \nu+\alpha +m+1 ) (2(k-m)+\alpha) d_{m,k}  - k ( \nu+\alpha +m+1 ) d_{m,k-1} -  ( \nu+\alpha +m+1 ) (k+\nu+\alpha+1) d_{m,k+1}$$

$$ - m  (\nu+\alpha +m+1 ) (m+\nu) (2k+1+\nu+\alpha ) d_{m-1,k} + m k (\nu+\alpha +m+1 ) (m+\nu) d_{m-1,k-1} $$

$$+ m (\nu+\alpha +m+1 ) (m+\nu) (k+\nu+\alpha+1) d_{m-1,k+1} $$

$$-  m(m+1) d_{m,k}  - (m+1) m (2k+1+\nu+\alpha ) (m+\nu) d_{m-1,k} +   m k (m+1)   (m+\nu) d_{m-1,k-1} $$

$$+  m (m+1) (m+\nu) (k+\nu+\alpha+1) d_{m-1,k+1} + (m+1)  (2k+1+\nu+\alpha ) d_{m,k} -  k (m+1) d_{m,k-1} $$

$$- (m+1) (k+\nu+\alpha+1)  d_{m,k+1} +   {k(k-1) (2m+1+\nu) \over  k+\nu+\alpha}  \left[ d_{m,k-2}  -   d_{m,k-1} \right] $$

$$ -  { k(k-1) (2k-3+\nu+\alpha) \over  k+\nu+\alpha} d_{m,k-2} +  { k(k-1)(k-2)  \over  k+\nu+\alpha} d_{m,k-3} + { k(k-1) (k+\nu+\alpha-1) \over  k+\nu+\alpha} d_{m,k-1} $$

$$+  {k (2m+\nu+k) (2k-1+\nu+\alpha) \over  k+\nu+\alpha}  d_{m,k-1} - {k(k-1)  (2m+\nu+k)  \over  k+\nu+\alpha}  d_{m,k-2}- k (2m+\nu+k) d_{m,k} $$

$$+  {m k (k-1)  (m+\nu) (2k-3+\nu+\alpha) \over  k+\nu+\alpha}  d_{m-1,k-2} - {k(k-1)(k-2) m  (m+\nu)  \over  k+\nu+\alpha}  d_{m-1,k-3}$$

$$- {k m (k-1) (m+\nu)  (k+\nu+\alpha-1) \over  k+\nu+\alpha} d_{m-1,k-1} $$

$$-  {m k (k-1)  (m+\nu) (2k-1+\nu+\alpha) \over  k+\nu+\alpha}  d_{m-1,k-1} + {k(k-1)^2 m  (m+\nu)  \over  k+\nu+\alpha}  d_{m-1,k-2}$$

$$+  k m (k-1) (m+\nu)   d_{m-1,k} +  {m k  (m+\nu) (2k-1+\nu+\alpha)^2 \over  k+\nu+\alpha}  d_{m-1,k-1}$$

$$ -  {m k(k-1)  (m+\nu) (2k-1+\nu+\alpha)\over  k+\nu+\alpha}  d_{m-1,k-2} - m k (m+\nu) (2k-1+\nu+\alpha) d_{m-1,k} $$

$$- {k(k-1) m  (m+\nu) (2k-3+\nu+\alpha) \over  k+\nu+\alpha}  d_{m-1,k-2}   +  {k(k-1)(k-2) m  (m+\nu)  \over  k+\nu+\alpha}  d_{m-1,k-3} $$

$$+  {m k (k-1)  (m+\nu) (k-1+\nu+\alpha) \over  k+\nu+\alpha}  d_{m-1,k-1}  $$

$$-  m k (m+\nu) (2k+1+\nu+\alpha) d_{m-1,k}   +  k^2  m  (m+\nu)   d_{m-1,k-1} $$

$$+  m k (m+\nu) (k+1+\nu+\alpha)  d_{m-1,k+1}  $$

$$-  { k  (2k-1+\nu+\alpha)^2 \over  k+\nu+\alpha}  d_{m,k-1} +  { k(k-1)  (2k-1+\nu+\alpha)\over  k+\nu+\alpha}  d_{m,k-2}
+  k (2k-1+\nu+\alpha) d_{m,k} $$

$$+ {k(k-1)  (2k-3+\nu+\alpha) \over  k+\nu+\alpha}  d_{m,k-2}   -  {k(k-1)(k-2)  \over  k+\nu+\alpha}  d_{m,k-3} -  { k (k-1)  (k-1+\nu+\alpha) \over  k+\nu+\alpha}  d_{m,k-1}  $$

$$+  k  (2k+1+\nu+\alpha) d_{m,k}   -  k^2 d_{m,k-1} -   k (k+1+\nu+\alpha)  d_{m,k+1} . $$
Hence after simplification we get  (3.20).

\end{proof}

On the other hand, taking into account orthogonality conditions (3.15), we have 

$$Q_{2n}(x)= \sum_{j=0}^{2n} c_{n,j}    L_j^{\nu+\alpha} (x) =  \sum_{j=0}^{n} c_{n,2j}  \   L_{2j}^{\nu+\alpha} (x)+  \sum_{j=0}^{n-1} c_{n,2j+1} \   L_{2j+1}^{\nu+\alpha} (x),$$
and via the uniqueness of the expansion of the associated polynomial $Q_{2n}$ by Laguerre polynomials we find

$$ c_{n,2j} = \sum_{m= j}^n  \ a_{n,m}  \ d_{m, 2j},\quad  c_{n,2j+1} = \sum_{m= j+1}^{n}   \ a_{n,m}  \ d_{m, 2j+1}.\eqno(3.22)$$  
We observe via (3.18)  that $d_{m, 2j} \neq 0, \  m= j,\dots, n, \  d_{m, 2j+1} \neq 0,  m=j+1,\dots, n.$ But from (3.15) we get for $n \in \mathbb{N}$

$$ c_{2n, 2j} =  0,\quad j=0,1,\dots  n-1 ;\quad  c_{2n,2j+1} = 0,\  j=0,1,\dots  n-1,$$ 

$$ c_{2n+1, 2j} =  0,\quad j=0,1,\dots  n ;\quad c_{2n+1,2j+1} = 0,\  j=0,1,\dots  n-1.$$ 
Consequently, equalities (3.22) represent for the polynomial sequence $\left(P_{2n}^{\nu,\alpha}\right)_{ n\in \mathbb{N}_0}\  \left(\left(P_{2n+1}^{\nu,\alpha}\right)_{ n\in \mathbb{N}_0} \right)$ linear homogeneous systems of $2n\  (2n+1)$ equations with  $2n+1\  (2(n+1) )$ unknowns.  However, if we assume that the free coefficient $a_{2n,0} \  ( a_{2n+1,0} ) $ is known, we come out with  linear non-homogeneous systems of  $2n\  (2n+1)$ equations with  $2n\  (2n+1 )$ unknowns.  It can be solved uniquely by Cramer's rule with nonzero determinant. In fact, we have the following non-homogeneous systems of $2n,\ 2n+1$ linear equations to determine the sequences  $\left(P_{2n}^{\nu,\alpha}\right)_{ n\in \mathbb{N}_0},\  \left(P_{2n+1}^{\nu,\alpha}\right)_{ n\in \mathbb{N}_0}  $, respectively,

$$\begin{pmatrix}

d_{1, 0}  & d_{2, 0}&  \dots&  \dots&  d_{n-1,0}&   d_{n,0} & \dots& d_{2n, 0} \\

d_{1, 1}&  d_{2, 1} &   \dots&    \dots&   \dots&    \dots &\dots& d_{2n, 1} \\
 
d_{1, 2}&  d_{2, 2} &   \dots&   \dots&   \dots&    \dots &\dots& d_{2n, 2} \\
 
0    &  d_{2, 3} &   \dots&   \dots& \dots&    \dots &  \dots& d_{2n, 3}\\
 
 \vdots    &  d_{2, 4} &   \dots&  \dots&   \dots&    \dots & \dots& d_{2n, 4}\\

 \vdots  &  0 &  d_{3, 5} &  \dots&   \dots&   \dots  &  \dots& d_{2n, 5}\\
 
 \vdots  &   \vdots  &  d_{3, 6} & \dots&   \dots&    \dots  & \dots&  d_{2n, 6}\\

   \vdots  &   \vdots  &  0 & \ddots&  \vdots&    \vdots  &  \ddots& \vdots\\
 
 \vdots&  \ddots & \ddots&  0 &   d_{n-1, 2n-3}  &  d_{n, 2n-3}&  \dots& d_{2n, 2n-3}\\    
 
  \vdots &  \ddots & \ddots&  0 &   d_{n-1, 2(n-1)}  &  d_{n, 2(n-1)}&   \dots&  d_{2n, 2(n-1)}\\  
  
  0&  \dots & \dots&   0 &  0 &  d_{n, 2n-1}  &  \dots&    d_{2n, 2n-1}\\  

\end{pmatrix} 
\begin{pmatrix}

a_{2n,1} \\
 
 a_{2n,2} \\

 \vdots  \\
 
  \vdots  \\
  
      \vdots  \\
    
     \vdots  \\
     
      \vdots  \\
      
       \vdots  \\
       
       a_{2n,2n-1}  \\
        
        a_{2n,2n}

\end{pmatrix}=  
\begin{pmatrix}

- a_{2n,0} \\
 
 0 \\

 \vdots  \\
 
  \vdots  \\
  
      \vdots  \\
    
     \vdots  \\
     
      \vdots  \\
      
       \vdots  \\
       
      0  \\
        
      0

\end{pmatrix}
 ,\eqno(3.23)$$
\vspace{0,5cm}

$$\begin{pmatrix}

d_{1, 0}  & d_{2, 0}&  \dots&  \dots& d_{n,0} &     \dots& d_{2n+1, 0} \\

d_{1, 1}&  d_{2, 1} &   \dots&   \dots&   \dots&    \dots& d_{2n+1, 1} \\
 
d_{1, 2}&  d_{2, 2} &   \dots&   \dots&   \dots&    \dots& d_{2n+1, 2} \\
 
0    &  d_{2, 3} &   \dots&   \dots& \dots&      \dots& d_{2n+1, 3}\\
 
 \vdots    &  d_{2, 4} &   \dots&   \dots&  \dots&     \dots& d_{2n+1, 4}\\

 \vdots  &  0 &  d_{3, 5} &  \dots&  \dots&     \dots& d_{2n+1, 5}\\
 
 \vdots  &   \vdots  &  d_{3, 6} & \dots&  \dots&     \dots&  d_{2n+1, 6}\\

   \vdots  &   \ddots  &  \ddots & \vdots&  \vdots&      \ddots& \vdots\\
 
 0&  \dots &  \dots &  0 &   d_{n, 2n-1}  &   \dots& d_{2n+1, 2n-1}\\    
 
  0&  \dots & \dots&   0 &   d_{n, 2n}  &     \dots&  d_{2n+1, 2n}\\  
  
 \end{pmatrix} 
\begin{pmatrix}

a_{2n+1,1} \\
 
  \vdots  \\
 
  \vdots  \\
  
      \vdots  \\
    
     \vdots  \\
     
      \vdots  \\
      
       \vdots  \\
       
       a_{2n+1,2n}  \\
        
        a_{2n+1,2n+1}

\end{pmatrix}=  
\begin{pmatrix}

- a_{2n+1,0} \\
 
 0 \\

 \vdots  \\
 
  \vdots  \\
  
      \vdots  \\
    
     \vdots  \\
     
      \vdots  \\
      
       \vdots  \\
                  
      0

\end{pmatrix}
 .\eqno(3.24)$$
Denoting by $D_{2n},\ D_{2n+1}$ the corresponding nonzero determinants of the systems (3.23), (3.24)

$$ D_{2n}= \begin{vmatrix}  

d_{1, 0}  & d_{2, 0}&  \dots&  \dots&  \dots&    \dots & \dots& d_{2n, 0} \\

d_{1, 1}&  d_{2, 1} &   \dots&    \dots&   \dots&    \dots &\dots& d_{2n, 1} \\
 
d_{1, 2}&  d_{2, 2} &   \dots&   \dots&   \dots&    \dots &\dots& d_{2n, 2} \\
 
0    &  d_{2, 3} &   \dots&   \dots& \dots&    \dots &  \dots& d_{2n, 3}\\
 
 \vdots    &  d_{2, 4} &   \dots&  \dots&   \dots&    \dots & \dots& d_{2n, 4}\\

 \vdots  &  0 &  d_{3, 5} &  \dots&   \dots&   \dots  &  \dots& d_{2n, 5}\\
 
 \vdots  &   \vdots  &  d_{3, 6} & \dots&   \dots&    \dots  & \dots&  d_{2n, 6}\\

   \vdots  &   \vdots  &  \vdots & \vdots&  \vdots&    \vdots  &  \ddots& \vdots\\
 
 0&  \dots & \dots&  0 &   d_{n-1, 2n-3}  &  \dots&  \dots& d_{2n, 2n-3}\\    
 
  0&  \dots & \dots&  0 &   d_{n-1, 2(n-1)}  &  \dots&   \dots&  d_{2n, 2(n-1)}\\  
  
  0&  \dots & \dots&   0 &  0 &  d_{n, 2n-1}  &  \dots&    d_{2n, 2n-1}\\  

\end{vmatrix},\eqno(3.25)$$

$$D_{2n+1}= \begin{vmatrix}

d_{1, 0}  & d_{2, 0}&  \dots&  \dots& \dots&     \dots& d_{2n+1, 0} \\

d_{1, 1}&  d_{2, 1} &   \dots&   \dots&   \dots&    \dots& d_{2n+1, 1} \\
 
d_{1, 2}&  d_{2, 2} &   \dots&   \dots&   \dots&    \dots& d_{2n+1, 2} \\
 
0    &  d_{2, 3} &   \dots&   \dots& \dots&      \dots& d_{2n+1, 3}\\
 
 \vdots    &  d_{2, 4} &   \dots&   \dots&  \dots&     \dots& d_{2n+1, 4}\\

 \vdots  &  0 &  d_{3, 5} &  \dots&  \dots&     \dots& d_{2n+1, 5}\\
 
 \vdots  &   \vdots  &  d_{3, 6} & \dots&  \dots&     \dots&  d_{2n+1, 6}\\

   \vdots  &   \vdots  &  \vdots & \vdots&  \vdots&      \ddots& \vdots\\
 
 0&  \dots &  \dots &  0 &   d_{n, 2n-1}  &   \dots& d_{2n+1, 2n-1}\\    
 
  0&  \dots & \dots&   0 &   d_{n, 2n}  &     \dots&  d_{2n+1, 2n}\\  
  
 \end{vmatrix},\eqno(3.26)$$ 
we apply Cramer's rule to get the expressions for the coefficients of the sequences  $\left(P_{2n}^{\nu,\alpha}\right)_{ n\in \mathbb{N}_0},\  \left(P_{2n+1}^{\nu,\alpha}\right)_{ n\in \mathbb{N}_0} $ in terms of the related free coefficients. Precisely, denoting by 

$$ D_{2n,1} = \begin{vmatrix}  

  d_{2, 1} &   \dots&    \dots&     \dots &   \dots &    \dots & d_{2n, 1} \\
 
 d_{2, 2} &   \dots&   \dots&   \dots&    \dots& \dots & d_{2n, 2} \\
 
  d_{2, 3} &   \dots&   \dots& \dots&      \dots&  \dots &d_{2n, 3}\\
 
  d_{2, 4} &   \dots&  \dots&       \dots & \dots&  \dots & d_{2n, 4}\\

   0 &  d_{3, 5} &  \dots&   \dots&    \dots&  \dots & d_{2n, 5}\\
 
    \vdots  &  d_{3, 6} & \dots&       \dots  & \dots&  \dots &  d_{2n, 6}\\

      \vdots  &  0 & \vdots&    \vdots  &  \ddots& \ddots & \vdots\\
 
 0 & \dots&  \dots&  \dots&   \dots&  \dots & d_{2n, 2n-3}\\    
 
  0 & \dots& \dots&    \dots&     \dots&  \dots & d_{2n,2(n-1)}\\  
  
  0 & \dots&   \dots&  \dots  &  \dots&   \dots &  d_{2n, 2n-1}\\  

\end{vmatrix},\eqno(3.27)$$

\vspace{2mm}

$$ D_{2n,k} = \begin{vmatrix}  

d_{1, 1}&  d_{2, 1} &   \dots&    \dots&    d_{k-1,1} &  d_{k+1,1} &    \dots & d_{2n, 1} \\
 
d_{1, 2}&  d_{2, 2} &   \dots&   \dots&   \dots&    \dots& \dots & d_{2n, 2} \\
 
0    &  d_{2, 3} &   \dots&   \dots& \dots&      \dots&  \dots &d_{2n, 3}\\
 
 \vdots    &  d_{2, 4} &   \dots&  \dots&       \dots & \dots&  \dots & d_{2n, 4}\\

 \vdots  &  0 &  d_{3, 5} &  \dots&   \dots&    \dots&  \dots & d_{2n, 5}\\
 
 \vdots  &   \vdots  &  d_{3, 6} & \dots&       \dots  & \dots&  \dots &  d_{2n, 6}\\

   \vdots  &   \vdots  &  \vdots & \vdots&    \vdots  &  \ddots& \ddots & \vdots\\
 
 0&  \dots & \dots&  \dots&  \dots&   \dots&  \dots & d_{2n, 2n-3}\\    
 
  0&  \dots & \dots& \dots&    \dots&     \dots&  \dots & d_{2n,2(n-1)}\\  
  
  0&  \dots & \dots&   \dots&  \dots  &  \dots&   \dots &  d_{2n, 2n-1}\\  

\end{vmatrix},\ k = 2,\dots, 2n-1,\eqno(3.28)$$

\vspace{2mm}

$$D_{2n,2n} = \begin{vmatrix}  

d_{1, 1}&  d_{2, 1} &   \dots&    \dots&    \dots&  \dots&    \dots & d_{2n-1, 1} \\
 
d_{1, 2}&  d_{2, 2} &   \dots&   \dots&   \dots&    \dots& \dots & d_{2n-1, 2} \\
 
0    &  d_{2, 3} &   \dots&   \dots& \dots&      \dots&  \dots &d_{2n-1, 3}\\
 
 \vdots    &  d_{2, 4} &   \dots&  \dots&       \dots & \dots&  \dots & d_{2n-1, 4}\\

 \vdots  &  0 &  d_{3, 5} &  \dots&   \dots&    \dots&  \dots & d_{2n-1, 5}\\
 
 \vdots  &   \vdots  &  d_{3, 6} & \dots&       \dots  & \dots&  \dots &  d_{2n-1, 6}\\

   \vdots  &   \vdots  &  \vdots & \vdots&    \vdots  &  \ddots& \ddots & \vdots\\
 
 0&  \dots & \dots&  \dots&  \dots&   \dots&  \dots & d_{2n-1, 2n-3}\\    
 
  0&  \dots & \dots& \dots&    \dots&     \dots&  \dots & d_{2n-1,2(n-1)}\\  
  
  0&  \dots & \dots&   \dots&  \dots  &  \dots&   \dots &  d_{2n-1, 2n-1}\\  

\end{vmatrix},\eqno(3.29)$$

\vspace{2mm}

$$ D_{2n+1,1}=  \begin{vmatrix}  

 d_{2, 1} &   \dots&    \dots&    \dots & \dots&    \dots & d_{2n+1, 1} \\
 
  d_{2, 2} &   \dots&   \dots&   \dots&    \dots& \dots & d_{2n+1, 2} \\
 
 d_{2, 3} &   \dots&   \dots& \dots&      \dots&  \dots &d_{2n+1, 3}\\
 
   d_{2, 4} &   \dots&  \dots&       \dots & \dots&  \dots & d_{2n+1, 4}\\

   0 &  d_{3, 5} &  \dots&   \dots&    \dots&  \dots & d_{2n+1, 5}\\
 
   \vdots  &  d_{3, 6} & \dots&       \dots  & \dots&  \dots &  d_{2n+1, 6}\\

    \vdots  &  \vdots & \vdots&    \vdots  &  \ddots& \ddots & \vdots\\
 
 0 & \dots&  \dots&  \dots&   \dots&  \dots & d_{2n+1, 2n-1}\\    
 
  0 & \dots& \dots&    \dots&     \dots&  \dots & d_{2n+1,2n}\\  
  
 \end{vmatrix},\eqno(3.30)$$

\vspace{2mm}

$$ D_{2n+1,k}=  \begin{vmatrix}  

d_{1, 1}&  d_{2, 1} &   \dots&    \dots&    d_{k-1,1} &  d_{k+1,1} &    \dots & d_{2n+1, 1} \\
 
d_{1, 2}&  d_{2, 2} &   \dots&   \dots&   \dots&    \dots& \dots & d_{2n+1, 2} \\
 
0    &  d_{2, 3} &   \dots&   \dots& \dots&      \dots&  \dots &d_{2n+1, 3}\\
 
 \vdots    &  d_{2, 4} &   \dots&  \dots&       \dots & \dots&  \dots & d_{2n+1, 4}\\

 \vdots  &  0 &  d_{3, 5} &  \dots&   \dots&    \dots&  \dots & d_{2n+1, 5}\\
 
 \vdots  &   \vdots  &  d_{3, 6} & \dots&       \dots  & \dots&  \dots &  d_{2n+1, 6}\\

   \vdots  &   \vdots  &  \vdots & \vdots&    \vdots  &  \ddots& \ddots & \vdots\\
 
 0&  \dots & \dots&  \dots&  \dots&   \dots&  \dots & d_{2n+1, 2n-1}\\    
 
  0&  \dots & \dots& \dots&    \dots&     \dots&  \dots & d_{2n+1,2n}\\  
  
 \end{vmatrix},\ k = 2,\dots, 2n,\eqno(3.31)$$
\vspace{2mm}

$$  D_{2n+1,2n+1}= \begin{vmatrix}  

d_{1, 1}&  d_{2, 1} &   \dots&    \dots&    \dots&  \dots&    \dots & d_{2n, 1} \\
 
d_{1, 2}&  d_{2, 2} &   \dots&   \dots&   \dots&    \dots& \dots & d_{2n, 2} \\
 
0    &  d_{2, 3} &   \dots&   \dots& \dots&      \dots&  \dots &d_{2n, 3}\\
 
 \vdots    &  d_{2, 4} &   \dots&  \dots&       \dots & \dots&  \dots & d_{2n, 4}\\

 \vdots  &  0 &  d_{3, 5} &  \dots&   \dots&    \dots&  \dots & d_{2n, 5}\\
 
 \vdots  &   \vdots  &  d_{3, 6} & \dots&       \dots  & \dots&  \dots &  d_{2n, 6}\\

   \vdots  &   \vdots  &  \vdots & \vdots&    \vdots  &  \ddots& \ddots & \vdots\\
 
 0&  \dots & \dots&  \dots&  \dots&   \dots&  \dots & d_{2n, 2n-1}\\    
 
  0&  \dots & \dots& \dots&    \dots&     \dots&  \dots & d_{2n,2n}\\  
  
 \end{vmatrix},\eqno(3.32)$$
we obtain the values  for coefficients of the sequences  $\left(P_{2n}^{\nu,\alpha}\right)_{ n\in \mathbb{N}_0}\  \left(P_{2n+1}^{\nu,\alpha}\right)_{ n\in \mathbb{N}_0} $, respectively, 

$$  a_{2n,k} = (-1)^{k}   a_{2n,0}\  { D_{2n,k}  \over D_{2n}},\quad   k = 1,\dots, 2n,\eqno(3.33)$$

$$  a_{2n+1,k} = (-1)^{k}   a_{2n+1,0}\  { D_{2n+1,k}  \over D_{2n+1}},\quad  k = 1,\dots, 2n+1.\eqno(3.34)$$
Moreover, returning to (3.5), we immediately obtain the values of the coefficients for the 3-term recurrence relation (3.4).  Indeed, we have

$$A_{2n+1} = - {a_{2n,0} \over a_{2n+1,0} } \  { D_{2n,2n} \  D_{2n+1}\over  D_{2n+1,2n+1}\   D_{2n} },\eqno(3.35)$$

$$A_{2n} = - {a_{2n-1,0} \over a_{2n,0} } \  { D_{2n-1,2n-1} \  D_{2n}\over  D_{2n,2n}\   D_{2n-1} },\eqno(3.36)$$

$$B_{2n} =  { D_{2n+1,2n}  \over D_{2n+1,2n+1}} -   { D_{2n,2n-1}  \over D_{2n,2n}},\eqno(3.37)$$  

$$B_{2n+1} =   { D_{2(n+1),2n+1}  \over D_{2(n+1),2(n+1)}} -  { D_{2n+1,2n}  \over D_{2n+1,2n+1}} .\eqno(3.38)$$  
In order to find free coefficients of the even and odd Prudnikov's sequences, we appeal to the identity (3.7) and values of the moments for $\rho_\nu$. Thus using (3.33),  we derive from (3.7) for the sequence $\left(P_{2n}^{\nu,\alpha}\right)_{ n\in \mathbb{N}_0}$

$${  D_{2n} \over  a_{2n,0}\  D_{2n,2n} } =   {a_{2n,0} \over  D_{2n}} \sum_{m=0}^{2n}  (-1)^{m}  D_{2n,m}  \Gamma(2n+ m+\alpha +\nu+1)\Gamma(2n+m+\alpha +1),\  D_{2n,0} \equiv D_{2n}.$$ 
Hence, taking into account the positive sign of the leading coefficient $a_{2n}$, we get the value of $a_{2n,0}$ in the form

$$ a_{2n,0} =  {  D_{2n} \over  \left[ D _{2n,2n}\right]^{1/2} } \left[ \sum_{m=0}^{2n}  (-1)^{m}  D_{2n,m}  \Gamma(2n+ m+\alpha +\nu+1)\Gamma(2n+m+\alpha +1)\right]^{-1/2},\ \  D_{2n,0} \equiv D_{2n}.\eqno(3.39) $$
Analogously, we obtain the value $a_{2n+1,0}$ for the odd sequence  $\left(P_{2n+1}^{\nu,\alpha}\right)_{n\in \mathbb{N}_0}$, namely,

$$ a_{2n+1,0} = -  {  D_{2n+1} \over  \left[ D _{2n+1,2n+1}\right]^{1/2} } \left[ \sum_{m=0}^{2n+1}  (-1)^{m}  D_{2n+1,m}  \Gamma(2(n+1) + m+\alpha +\nu)\Gamma(2(n+1)+m+\alpha )\right]^{-1/2},\eqno(3.40)$$
 where $  D_{2n+1,0} \equiv  D_{2n+1}. $  Leading coefficients for the Prudnikov sequences have the values, accordingly,

$$ a_{2n} =   \left[ D _{2n,2n}\right]^{1/2} \left[ \sum_{m=0}^{2n}  (-1)^{m}  D_{2n,m}  \Gamma(2n+ m+\alpha +\nu+1)\Gamma(2n+m+\alpha +1)\right]^{-1/2},\eqno(3.41) $$

$$ a_{2n+1} =  \left[ D _{2n+1,2n+1}\right]^{1/2} \left[ \sum_{m=0}^{2n+1}  (-1)^{m}  D_{2n+1,m}  \Gamma(2(n+1) + m+\alpha +\nu)\Gamma(2(n+1)+m+\alpha )\right]^{-1/2}.\eqno(3.42) $$
Thus we proved the following theorem.

{\bf Theorem 2.} {\it Let $\nu \ge 0,\ \alpha > -1,\ n \in \mathbb{N}_0.$ Prudnikov's sequences of orthogonal polynomials  $\left(P_{2n}^{\nu,\alpha}\right)_{ n\in \mathbb{N}_0},\\  \left(P_{2n+1}^{\nu,\alpha}\right)_{ n\in \mathbb{N}_0} $ have explicit values with coefficients calculated by formulas 

$$  a_{2n,k} = (-1)^{k}   a_{2n,0}\  { D_{2n,k}  \over D_{2n}},\quad   k = 1,\dots, 2n,$$

$$  a_{2n+1,k} = (-1)^{k}   a_{2n+1,0}\  { D_{2n+1,k}  \over D_{2n+1}},\quad  k = 1,\dots, 2n+1,$$
respectively, where  the determinants  $ D_{2n}, D_{2n+1}, D_{2n,k}, D_{2n+1,k}$ are defined by $(3.25)- (3.32)$ and free coefficients  $a_{2n,0}, a_{2n+1,0}$ by  $(3.37), (3.38)$. Moreover, the 3-term recurrence relation $(3.4)$ holds with coefficients $(3.35)- (3.38)$.}

{\bf Corollary 1}.  {\it Coefficients $(3.14)$ are calculated by  formulas

$$ c_{2n,2j} = {a_{2n,0} \over D_{2n}} \sum_{m= j}^{2n}  \  (-1)^m D_{2n,m}  \ d_{m, 2j},\ j= n,\dots, 2n,\eqno(3.43)$$

$$ c_{2n+1,2j} =  {a_{2n+1,0} \over D_{2n+1}} \sum_{m= j}^{2n+1}  \   (-1)^m D_{2n+1,m}   \ d_{m, 2j}, \ j= n+1,\dots, 2n+1, \eqno(3.44)$$

$$ c_{2n,2j+1} =    {a_{2n,0} \over D_{2n}}   \sum_{m= j+1}^{2n}  ( -1)^m D_{2n,m}  \ d_{m, 2j+1}, \ j= n,\dots, 2n -1, \eqno(3.45)$$  

$$ c_{2n+1,2j+1} =   {a_{2n+1,0} \over D_{2n+1}}  \sum_{m= j+1}^{2n+1}  \  (-1)^m D_{2n+1,m}   \ d_{m, 2j+1},\  j= n,\dots, 2n.\eqno(3.46)$$  
where values  $ D_{2n}, D_{2n+1}, D_{2n,k}, D_{2n+1,k}$ are defined by $(3.25)- (3.32)$ and free coefficients  $a_{2n,0}, a_{2n+1,0}$ by 
$(3.39), (3.40)$}.

Our goal now is to find an analog of the Rodrigues formula for Prudnikov's polynomials. To do this,  we recall the representation (2.23) of an arbitrary polynomial in terms of its associated polynomial and representations (2.28), (2.29), (3.15), (3.16)  to write the following equalities  for the sequence $\left(P_n^{\nu,\alpha}\right)_{n\in \mathbb{N}_0}$

$$ P_n^{\nu,\alpha} (x) = {x^{-\alpha} \over \rho_\nu(x) }   \sum_{j=n}^{2n} {c_{n,j} \over j! }  S_j^{\nu,\alpha} (x)  = {x^{-\alpha} \over \rho_\nu(x) }   \sum_{j=n}^{2n} {c_{n,j} \over j! } {d^j\over dx^j} \left[ x^{j+\alpha} \rho_\nu(x)\right] $$ 

$$ = {x^{-\alpha} \over \rho_\nu(x) }   \sum_{j=0}^{n} {c_{n,j+n} \over (j+n)! } {d^{j+n}\over dx^{j+n}} \left[ x^{j+n+\alpha} \rho_\nu(x)\right] .$$ 
Therefore for sequences $\left(P_{2n}^{\nu,\alpha}\right)_{ n\in \mathbb{N}_0},\  \left(P_{2n+1}^{\nu,\alpha}\right)_{ n\in \mathbb{N}_0} $  we have,  correspondingly, 

$$ P_{2n}^{\nu,\alpha} (x) =  {x^{-\alpha} \over \rho_\nu(x) }    \sum_{j=0}^{2n} {c_{2n,j+2n} \over (j+2n)! } S_{j+2n}^{\nu, \alpha} (x), \eqno(3.45)$$ 

$$ P_{2n+1}^{\nu,\alpha} (x) =  {x^{-\alpha} \over \rho_\nu(x) }   \sum_{j=0}^{2n+1} {c_{2n+1,j+2n+1} \over (j+2n+1)! } S_{j+2n+1}^{\nu, \alpha} (x). \eqno(3.48)$$
In the meantime,  the sums in (3.47), (3.48) can be treated as follows

$$ \sum_{j=0}^{2n} {c_{2n,j+2n} \over (j+2n)! } S_{j+2n}^{\nu, \alpha} (x) = \sum_{j=0}^{n} {c_{2n,2(j+n)} \over (2(j+n))! } S_{2(j+n)}^{\nu, \alpha} (x) +\sum_{j=0}^{n-1} {c_{2n, 2(j+n)+1} \over (2(j+n)+1)! } S_{2(j+n)+1}^{\nu, \alpha} (x),$$

$$ \sum_{j=0}^{2n+1} {c_{2n+1,j+2n+1} \over (j+2n+1)! } S_{j+2n+1}^{\nu, \alpha} (x) =   \sum_{j=0}^{n} {c_{2n+1, 2(j+n+1)} \over (2(j+n+1))! } S_{2(j+n+1)}^{\nu, \alpha} (x) + \sum_{j=0}^{n} {c_{2n+1, 2(j+n)+1} \over (2(j+n)+1)! } S_{2(j+n)+1}^{\nu, \alpha} (x).$$
Hence, employing the theory of multiple orthogonal polynomials associated with the scaled Macdonald functions and the related  Rodrigues formulas (see details in \cite{AsscheYakubov2000}, \cite{Cous}), we find the following expressions

$$S_{2(j+n)}^{\nu, \alpha} (x) = x^{\alpha} \left[ A^{\alpha}_{j+n, j+n-1}(x)  \rho_\nu(x)+  B^{\alpha}_{j+n, j+n-1} (x)\rho_{\nu+1} (x) \right],\eqno(3.49)$$

$$S_{2(j+n)+1}^{\nu, \alpha} (x) = x^{\alpha} \left[ A^{\alpha}_{j+n, j+n} (x) \rho_\nu(x)+  B^{\alpha}_{j+n, j+n} (x)\rho_{\nu+1} (x) \right],\eqno(3.50)$$
where $A$-polynomials in front of $\rho_\nu$ are of degree $j+n$ as well as $B$-polynomial in (3.50), while $B$-polynomial in (3.49) is  of degree $j+n-1$.   These polynomials are explicitly calculated in \cite{Cous}.  Therefore formulas (3.47), (3.48) become, respectively,

$$ P_{2n}^{\nu,\alpha} (x) =   \sum_{j=0}^{n}   {c_{2n,2(j+n)} \over (2(j+n))! }  A^{\alpha}_{j+n, j+n-1}(x)  +  \sum_{j=0}^{n-1} {c_{2n, 2(j+n)+1} \over (2(j+n)+1)! }  A^{\alpha}_{j+n, j+n} (x) $$

$$+ {\rho_{\nu+1} (x) \over \rho_\nu(x)} \left[ \sum_{j=0}^{n} {c_{2n,2(j+n)} \over (2(j+n))! }  B^{\alpha}_{j+n, j+n-1}(x)  +  \sum_{j=0}^{n-1} {c_{2n, 2(j+n)+1} \over (2(j+n)+1)! }  B^{\alpha}_{j+n, j+n} (x) \right], \eqno(3.51)$$

$$ P_{2n+1}^{\nu,\alpha} (x) = \sum_{j=0}^{n}  \left[ {c_{2n+1,2(j+n+1)} \over (2(j+n+1))! }  A^{\alpha}_{j+n+1, j+n}(x)  +  {c_{2n+1, 2(j+n)+1} \over (2(j+n)+1)! }  A^{\alpha}_{j+n, j+n} (x) \right] $$

$$ + {\rho_{\nu+1} (x) \over \rho_\nu(x)} \sum_{j=0}^{n}  \left[  {c_{2n+1,2(j+n+1)} \over (2(j+n+1))! }  B^{\alpha}_{j+n+1, j+n}(x)  +   {c_{2n+1, 2(j+n)+1} \over (2(j+n)+1)! }  B^{\alpha}_{j+n, j+n} (x) \right]. \eqno(3.52)$$ 
But Lemma 1 presumes immediately the following identities from (3.51), (3.52)

$$ P_{2n}^{\nu,\alpha} (x) =   \sum_{j=n}^{2n}   {c_{2n,2j} \over (2j)! }  A^{\alpha}_{j, j-1}(x)  +  \sum_{j=n}^{2n-1} {c_{2n, 2j+1} \over (2j+1)! }  A^{\alpha}_{j, j} (x) ,\eqno(3.53)$$

$$ P_{2n+1}^{\nu,\alpha} (x) = \sum_{j=n}^{2n}  \left[ {c_{2n+1,2(j+1)} \over (2(j+1))! }  A^{\alpha}_{j+1, j}(x)  +  {c_{2n+1, 2j+1} \over (2j+1)! }  A^{\alpha}_{j, j} (x) \right] ,\eqno(3.54)$$
giving explicit expressions of Prudnikov's polynomials in terms of the multiple orthogonal polynomials for the scaled Macdonald functions, and two more relations between multiple $B$-polynomials

$$ \sum_{j=n}^{2n} {c_{2n,2j} \over (2j)! }  B^{\alpha}_{j, j-1}(x)  +  \sum_{j=n}^{2n-1} {c_{2n, 2j+1} \over (2j+1)! }  B^{\alpha}_{j, j} (x) \equiv 0, \eqno(3.55)$$ 

$$\sum_{j=n}^{2n}  \left[  {c_{2n+1,2(j+1)} \over (2(j+1))! }  B^{\alpha}_{j+1, j}(x)  +   {c_{2n+1, 2j+1} \over (2j+1)! }  B^{\alpha}_{j, j} (x) \right] \equiv 0. \eqno(3.56)$$ 
On the other hand, 

$$ P_n^{\nu,\alpha} (x) = {x^{-\alpha} \over \rho_\nu(x) }   \sum_{j=n}^{2n} {c_{n,j} \over j! }  S_j^{\nu,\alpha} (x)  = {x^{-\alpha} \over \rho_\nu(x) }  {d^n\over dx^n}  \sum_{j=0}^{n} {c_{n,j+n} \over (j+n)! } S_j^{\nu, n+\alpha} (x). \eqno(3.57)$$ 
Hence, recalling integral representations (2.3),  (2.29) and the explicit formula for the associated Laguerre polynomials \cite{Chi}, we obtain
$$S^{\nu,n+\alpha}_j(x) =   x^{n+\alpha} j! \sum_{k=0}^j {(-1)^k\over k!} \binom{j+n+\nu+\alpha}{j-k} \int_0^\infty  e^{-t-x/t} t^{\nu+k-1} dt$$ 

$$=    x^{n+\alpha} j! \sum_{k=0}^j {(-1)^k\over k!} \binom{j+n+\nu+\alpha}{j-k} \rho_{\nu+k}(x).\eqno(3.58)$$
The problem now is to express $\rho_{\nu + k},\ k \in \mathbb{N}_0$ in terms of $\rho_{\nu}$ and $\rho_{\nu+1}$.   To do this, we use the Mellin-Barnes representation (1.4) and the definition (2.7) of the Pochhammer symbol  to derive
$$\rho_{\nu + k}(x) = {1\over 2\pi i } \int_{\gamma -i\infty}^{\gamma+i\infty} \Gamma(s+\nu+k)\Gamma(s) x^{-s} ds
=  {1\over 2\pi i } \int_{\gamma -i\infty}^{\gamma+i\infty} (s+\nu)_k \Gamma(s+\nu)\Gamma(s) x^{-s} ds $$
$$=  {(-1)^k x^{\nu+k} \over 2\pi i } {d^k\over dx^k}   \int_{\gamma -i\infty}^{\gamma+i\infty}  \Gamma(s+\nu)\Gamma(s) x^{-s-\nu} ds =    (-1)^k x^{\nu+k}  {d^k\over dx^k} \left[  x^{-\nu}  \rho_{\nu}(x) \right].$$
Then, employing  the Leibniz formula and (2.16), we find

$$ \rho_{\nu + k}(x) =  \sum_{m=0}^k \binom{k}{m}   (\nu)_{k-m} \  x^{m}  \rho_{\nu-m}(x). \eqno(3.59)$$
Meanwhile, employing the identity from \cite{Cous}  for the scaled Macdonald functions, precisely,

$$ x^{ m}  \rho_{\nu-m}(x) = x^{ m/2} r_m(2\sqrt x; \nu ) \rho_{\nu}(x) +  x^{ (m-1)/2} r_{m-1} (2\sqrt x; \nu -1) \rho_{\nu+1}(x), \quad  m \in {\mathbb N}_0,\eqno(3.60)$$
where $r_{-1}(z;\nu)=0$, 

$$ x^{ m/2} r_m(2\sqrt x; \nu ) = (-1)^m \sum_{i=0}^{[m/2]}   (\nu+i-m+1)_{m-2i} (m-2i+1)_i  {x^i\over i!},\eqno(3.61)$$
formula (3.59) takes the final expression 

$$  \rho_{\nu + k}(x) =    \rho_{\nu}(x) \sum_{m=0}^k \sum_{i=0}^{[m/2]}   (-1)^{m} \  (\nu+i-m+1)_{m-2i} \ (m-2i+1)_i \    (\nu)_{k-m}\ \binom{k}{m}   { x^i \over i!} $$

$$+   \rho_{\nu+1}(x) \sum_{m=0}^{k-1} \sum_{i=0}^{[m/2]}  (-1)^{m} \  (\nu+i-m)_{m-2i} \ (m-2i+1)_i \    (\nu)_{k-m-1}\ \binom{k}{m+1}   {x^i \over i!} .\eqno(3.62)$$
Substituting the right-hand side of the equality (3.62) into (3.58),  we get finally

$$S^{\nu,n+\alpha}_j(x) =   x^{n+\alpha} j! \left[   \rho_{\nu}(x) \sum_{k=0}^j \sum_{m=0}^k \sum_{i=0}^{[m/2]}  {(-1)^{k+m} \over k!} \binom{j+n+\nu+\alpha}{j-k}  (\nu+i-m+1)_{m-2i} \ (m-2i+1)_i \    (\nu)_{k-m}\ \binom{k}{m}   { x^i \over i!} \right.$$

$$\left. +   \rho_{\nu+1}(x) \sum_{k=0}^j \sum_{m=0}^{k-1} \sum_{i=0}^{[m/2]}   {(-1)^{k+m} \over k!} \binom{j+n+\nu+\alpha}{j-k} \  (\nu+i-m)_{m-2i} \ (m-2i+1)_i \    (\nu)_{k-m-1}\ \binom{k}{m+1}   {x^i \over i!} \right].\eqno(3.63)$$
Thus, returning to (3.57), we end up with the so-called Rodrigues type formula for the Prudnikov orthogonal polynomials $P_n^{\nu,\alpha}$

 $$ P_n^{\nu,\alpha} (x) = {x^{-\alpha} \over \rho_\nu(x) }  {d^n\over dx^n}  \left[  x^{n+\alpha}   \left[ \rho_{\nu}(x) \sum_{j=0}^{n}  \sum_{k=0}^j \sum_{m=0}^k \sum_{i=0}^{[m/2]} { (-1)^{k+m}   \over (j+n)! } \  c_{n, j+n}     (n+\nu+\alpha+k+1)_{j-k}  (\nu+i-m+1)_{m-2i}  \  \right. \right. $$
 
 $$ \left. \left. \times  \ (m-2i+1)_i (\nu)_{k-m}\  \binom{j}{k} \ \binom{k}{m}   { x^i \over i!}  +   \rho_{\nu+1}(x) \sum_{j=0}^{n} \sum_{k=0}^j \sum_{m=0}^{k-1} \sum_{i=0}^{[m/2]}  { (-1)^{k+m}   \over (j+n)! } \  c_{n, j+n}     (n+\nu+\alpha+k+1)_{j-k} \right. \right.$$
 
 $$\left.\left. \times \  (\nu+i-m)_{m-2i} \ (m-2i+1)_i \    (\nu)_{k-m-1}\  \binom{j}{k} \binom{k}{m+1}   {x^i \over i!} \right] \right]. \eqno(3.64)$$

{\bf Theorem 3.}  {\it Prudnikov's orthogonal polynomials $P_n^{\nu,\alpha} $ can be obtained from the Rodrigues type formula $(3.64)$, where coefficients $c_{n, j+n} $ are calculated in Corollary $1$. Moreover, Prudnikov's sequences $\left(P_{2n}^{\nu,\alpha}\right)_{ n\in \mathbb{N}_0}, \  \left(P_{2n+1}^{\nu,\alpha}\right)_{ n\in \mathbb{N}_0} $ are expressed in terms of multiple orthogonal polynomials related to the scaled Macdonald functions by equalities $(3.53), (3.54)$, respectively, where the polynomials $A_{j,j-1}^\alpha,  A_{j,j}^\alpha$ are calculated explicitly in \cite{Cous} by formulas}

$$ A_{j,j-1}^\alpha (x) = (\alpha+1)_{2j} \sum_{ m=0}^j  \binom{2j}{2m} {x^m\over  (\alpha+1)_{2m}} \ {}_3F_2 \left( - 2(j-m),\  m- \nu,\ m+1;\ 2m+1+\alpha,\ 2m+1;\ 1\right) ,$$

$$ A_{j,j}^\alpha (x) = (\alpha+1)_{2j+1} \sum_{ m=0}^j  \binom{2j+1}{2m} {x^m\over  (\alpha+1)_{2m}} \ {}_3F_2 \left( - 2(j-m)-1,\  m- \nu,\ m+1;\ 2m+1+\alpha,\ 2m+1;\ 1\right).$$

Further, the generating function for  polynomials $P_n^{\nu,\alpha}$ can be defined as usual by the equality

$$G(x,z) = \sum_{n=0}^\infty P_n^{\nu,\alpha} (x) {z^n\over n!} ,\quad  x >0,\  z \in \mathbb{C},\eqno(3.65)$$
where $|z| < h_x$ and $h_x >0$ is a  convergence radius of the power series. Then returning to (3.57)  and employing (2.28), we have from (3.65)

$$  G(x,z) =  {x^{-\alpha} \over \rho_\nu(x) } \sum_{n=0}^\infty  {z^n\over n!}  \sum_{j=n}^{2n} {c_{n,j} \over j! }   {d^{j}\over dx^{j}} \left[  x^{j+\alpha}  \rho_\nu(x) \right]=  {1\over \rho_\nu(x) } \sum_{n=0}^\infty  {z^n\over n!}  \sum_{j=n}^{2n} \sum_{k=0}^j  (-1)^k {c_{n,j} \over k! }   \binom{j+\alpha }{j-k} x^k \rho_{\nu-k}(x) . $$
Hence substituting the value of $x^k\rho_{\nu-k}(x)$ by formula (3.60), we get, finally, the expression for the generating function for the Prudnikov sequence $\left(P_n^{\nu,\alpha} \right)_{n\in\mathbb{N}_0}$, namely,

$$  G(x,z) =   \sum_{n=0}^\infty  \sum_{j=n}^{2n} \sum_{k=0}^j { (-1)^k \ c_{n,j} \over n!\ k! }   \binom{j+\alpha }{j-k} x^{k/2}  r_k(2\sqrt x; \nu ) z^n $$

$$+  {\rho_{\nu+1}(x) \over \rho_\nu(x) } \sum_{n=0}^\infty \sum_{j=n}^{2n} \sum_{k=0}^j  { (-1)^k c_{n,j} \over n!\ k! }   \binom{j+\alpha }{j-k} x^{(k-1)/2} r_{k-1}(2\sqrt x; \nu-1 ) z^n,\eqno(3.66) $$
where $c_{n,j}$ are defined in Corollary 1.

\section{Orthogonal polynomials with ultra-exponential weights}

In this section we will consider a sequence of  polynomials $\left(Q_n^{\nu,k}\right)_{n\in \mathbb{N}_0}$, which is  orthogonal with respect to the weight function (2.1) $x^\alpha \rho_{\nu,k}(x)$

$$\int_0^\infty  Q_n^{\nu,\alpha}(x) Q_m^{\nu,\alpha}(x)  \rho_{\nu,k}(x) x^\alpha dx = \delta_{m,n},\ \nu \ge0,\ \alpha > -1.\eqno(4.1)$$
The function $\rho_{\nu,k}$ satisfies some interesting properties.  In fact, recalling the Mellin-Barnes integral representation (2.1), we write

$$\rho_{\nu+1,k}(x)=  \frac{1}{2\pi i} \int_{\gamma-i\infty}^{\gamma+i\infty} \Gamma(\nu+1+s) \left[ \Gamma (s) \right]^k  x^{-s} ds$$

$$=  \frac{\nu}{2\pi i} \int_{\gamma-i\infty}^{\gamma+i\infty} \Gamma(\nu+s) \left[ \Gamma (s) \right]^k  x^{-s} ds+   \frac{1}{2\pi i} \int_{\gamma-i\infty}^{\gamma+i\infty} \Gamma(\nu+s)  s \left[ \Gamma (s) \right]^k  x^{-s} ds$$

$$=  \nu \rho_{\nu,k}(x) - x D \rho_{\nu,k}(x).$$
Hence, as in (2.18)

$$\rho_{\nu+1,k}(x)= (\nu- xD) \rho_{\nu,k}(x).\eqno(4.2)$$
Further,

$$ D \left( x D\right)^{k-1} \left( x^{\nu+1} D \left( x^{-\nu} \rho_{\nu,k}(x) \right) \right) = \frac{(-1)^{k+1}}{2\pi i} \int_{\gamma-i\infty}^{\gamma+i\infty} \Gamma(\nu+s+1)  s^k \left[ \Gamma (s) \right]^k  x^{-s-1} ds =   (-1)^{k+1} \rho_{\nu,k}(x).$$
Thus we derive the following $k+1$-th order differential equation for the function $\rho_{\nu,k}$

$$   (-1)^{k+1} D \left( x D\right)^{k-1} \left( x^{\nu+1} D \left( x^{-\nu} \rho_{\nu,k}(x) \right) \right) = \rho_{\nu,k} (x),\ k \in \mathbb{N},\quad  D \equiv {d\over dx}.\eqno(4.3)$$
From the Parseval equality (1.10)  it follows the integral recurrence relation for functions  $\rho_{\nu,k}$. Precisely, we obtain

$$\rho_{\nu,k+1}(x) = \int_0^\infty e^{-x/t} \rho_{\nu,k} (t) {dt\over t},\quad k \in \mathbb{N}_0.\eqno(4.4)$$
An analog of the integral representation (2.11) for   $\rho_{\nu,k}$ can be deduced in the following manner. In fact, the Mellin-Barnes integral for the associated Laguerre polynomials (cf. \cite{PrudnikovMarichev}, Vol. III, relation (8.4.33.3)) 

$$n!\  e^{-x}L_n^\nu(x)= {1\over 2\pi i}  \int_{\gamma-i\infty}^{\gamma+i\infty} \Gamma(s)\frac{ \Gamma(1+n+\nu-s)}{ \Gamma(1+\nu-s)} x^{-s} ds,\eqno(4.5)$$
integral (2.1) with the Parseval identity (1.10) and the reflection formula for the gamma-function imply the equality for $k \in \mathbb{N}$

$$x^n \rho_{\nu,k}(x) =  {1\over 2\pi i}  \int_{\gamma-i\infty}^{\gamma+i\infty} \Gamma(s+\nu+n) \left[  \Gamma(s+n)\right]^k x^{-s} ds$$

$$=  {1\over 2\pi i}  \int_{\gamma-i\infty}^{\gamma+i\infty} \Gamma(s+\nu+n) \Gamma(s+n) \Gamma(1-s-n) {\left[  \Gamma(s+n)\right]^{k-1}\over \Gamma(1-s-n)} x^{-s} ds$$

$$=  {(-1)^n \over 2\pi i}  \int_{\gamma-i\infty}^{\gamma+i\infty} \Gamma(s+\nu+n) \Gamma(s) \Gamma(1-s) {\left[  \Gamma(s+n)\right]^{k-1}\over \Gamma(1-s-n)} x^{-s} ds$$

$$=  (-1)^n n! \int_0^\infty  t^{\nu+n-1} e^{-t} L_n^\nu(t) \varphi_n \left({x\over t} \right) dt,$$
where

$$\varphi_{n,k} (x)=   {1\over 2\pi i}  \int_{\gamma-i\infty}^{\gamma+i\infty}  \Gamma(s) \left[  \Gamma(s+n)\right]^{k-1} x^{-s} ds.\eqno(4.6)$$
Therefore we obtain the integral representation

$$ x^n \rho_{\nu,k}(x) = (-1)^n n! \int_0^\infty  t^{\nu+n-1} e^{-t} L_n^\nu(t) \varphi_{n,k}  \left({x\over t} \right) dt.\eqno(4.7)$$
Differentiating (4.6) $n$ times by $x$, where the differentiation under the integral sign is possible due to the absolute and uniform convergence, we take into  the reduction formula for the gamma-function and (2.1) to obtain 

$$ {d^n\over dx^n } \varphi_{n,k} (x)=   {(-1)^n \over 2\pi i}  \int_{\gamma-i\infty}^{\gamma+i\infty} (s)_n \Gamma(s) \left[  \Gamma(s+n)\right]^{k-1} x^{-s-n} ds =  {(-1)^n \over 2\pi i}  \int_{\gamma-i\infty}^{\gamma+i\infty} \left[  \Gamma(s+n)\right]^{k} x^{-s-n} ds$$

$$= (-1)^n  \rho_{0,k-1}(x).\eqno(4.8)$$
Consequently, after  differentiation both sides of (4.7) $n$ times   we find an analog of the representation (2.29), namely,

$$   {d^n\over dx^n } \left[ x^n \rho_{\nu,k}(x)\right]  =  n! \int_0^\infty  t^{\nu-1} e^{-t} L_n^\nu(t) \rho_{0,k-1} \left({x\over t} \right) dt.\eqno(4.9)$$
Now, returning to (4.1), we substitute the  function  $\rho_{\nu,k}$ by the integral (4.4) and interchange the order of integration by Fubini theorem.   Then employing again the Viskov type identities (2.6) for the differential operator $\theta$, we derive for $k \in \mathbb{N}$

$$  \delta_{m,n} =  \int_0^\infty  Q_n^{\nu,\alpha}(x) Q_m^{\nu,\alpha}(x)  \rho_{\nu,k}(x) x^\alpha dx =   \int_0^\infty  \rho_{\nu,k-1} (t) {1\over t} \int_0^\infty  e^{-x/t}  Q_n^{\nu,\alpha}(x) Q_m^{\nu,\alpha}(x)  x^\alpha dx dt$$

$$= \Gamma(1+\alpha)  \int_0^\infty  \rho_{\nu,k-1} (t)  Q_n^{\nu,\alpha}(\theta ) Q_m^{\nu,\alpha}(\theta) \{ t^\alpha \} dt .$$  
Hence it leads to

{\bf Theorem 4.} {\it Let $k \in \mathbb{N},\ \nu \ge 0, \alpha > -1$.  The orthogonality $(4.1)$ for the sequence of polynomials  $\left(Q_n^{\nu,\alpha}\right)_{n\in \mathbb{N}_0} $ with the weight $x^\alpha \rho_{\nu,k}(x)$ is  the composition orthogonality of the same sequence with respect to the weight $ \rho_{\nu,k-1} $,  namely

$$ \int_0^\infty  \rho_{\nu,k-1} (t)  Q_n^{\nu,\alpha}(\theta ) Q_m^{\nu,\alpha}(\theta) \{ t^\alpha \} dt = {\delta_{m,n} \over \Gamma(1+\alpha)},\ m,n \in \mathbb{N}_0.\eqno(4.10)$$
In particular,  for $k=2$ this sequence is compositionally orthogonal in the sense of Prudnikov. }

Further, up to a normalization constant equality (4.1) is equivalent to the following $n$ conditions

$$\int_0^\infty  Q_n^{\nu,\alpha}(x) \rho_{\nu,k}(x) x^{\alpha+m} dx = 0,\ m= 0, 1,\dots, n-1,\ n \in \mathbb{N}.\eqno(4.11)$$
Hence the composition orthogonality (4.10) implies with the integration by parts and properties of the operator $\theta$

$$ \int_0^\infty   \theta^m \left\{ \rho_{\nu,k-1} (t) \right\}   Q_n^{\nu,\alpha} ( \theta ) \left\{ t^{\alpha}\right\} dt =  0,\ m= 0, 1,\dots, n-1,\  k, n \in \mathbb{N}.\eqno(4.12)$$
Writing $Q_n^{\nu,\alpha}$ in the explicit form 

$$ Q_n^{\nu,\alpha} (x) = \sum_{j=0}^n a_{n,j} x^j\equiv  Q_n^{\nu,\alpha,0} (x),$$
we have

$$ Q_n^{\nu,\alpha} ( \theta ) \left\{ t^{\alpha}\right\} =  \sum_{j=0}^n a_{n,j} \theta^j \left\{ t^{\alpha}\right\} = t^\alpha  \sum_{j=0}^n a_{n,j} (1+\alpha)_j \ t^j=  t^\alpha Q_n^{\nu,\alpha,1} ( t) ,$$
where

$$Q_n^{\nu,\alpha,1} ( t) = {1\over \Gamma(1+\alpha)}  \sum_{j=0}^n a_{n,j} \Gamma(1+\alpha+j) \ t^j.\eqno(4.13)$$

On the other hand, employing (4.4) for $k\ge 2$ and observing owing to the Viskov type identities (2.6) that  $(  \theta_t\equiv tDt,\ \beta_y\equiv DyD)$

$$\theta^m_t\left\{ e^{-ty}\right\} = (-1)^m \beta^m_y \left\{ e^{-ty}\right\}, \quad m \in \mathbb{N}_0,\eqno(4.14) $$
we deduce, integrating by parts

$$  \theta_t^m \left\{ \rho_{\nu,k-1} (t) \right\}  = \theta_t^m \left\{ \int_0^\infty e^{-t y} \rho_{\nu,k-2} \left({1\over y}\right) {dy\over y} \right\} $$

$$= (-1)^m  \int_0^\infty \beta^m_y \left\{ e^{-ty}\right\} \  \rho_{\nu,k-2} \left({1\over y}\right) {dy\over y} $$

$$= (-1)^m  \int_0^\infty  e^{-ty} \beta^m_y \left\{  \rho_{\nu,k-2} \left({1\over y}\right)  {1\over y} \right\} dy,$$
where the differentiation under integral sign is allowed via the absolute and uniform convergence. Thus, returning to (4.12), we  plug in   the latter expressions and change the order of integration by Fubini's theorem to write it in the form  

$$ \int_0^\infty     Q_n^{\nu,\alpha,2} \left( {1\over y}\right) \beta^m_y \left\{  \rho_{\nu,k-2} \left({1\over y}\right)  {1\over y} \right\}  y^{-\alpha-1} dy  =  0,\ m= 0, 1,\dots, n-1,\   n \in \mathbb{N},\eqno(4.15)$$
where

$$ Q_n^{\nu,\alpha,2} (x) =  {1\over \Gamma(1+\alpha)} \sum_{j=0}^n a_{n,j} [\Gamma (1+\alpha+j) ]^2  \ x^j.\eqno(4.16)$$
Meanwhile, recalling (2.1), we get 

$$\beta^m_y \left\{  \rho_{\nu,k-2} \left({1\over y}\right)  {1\over y} \right\}  = \left(D^m y^m D^m\right) \left\{ {1\over 2\pi i}  \int_{\gamma-i\infty}^{\gamma+i\infty}  \Gamma(s+\nu) \left[  \Gamma(s)\right]^{k-2} y^{s-1} ds\right\} $$
 
$$= {1\over 2\pi i}  \int_{\gamma-i\infty}^{\gamma+i\infty} [(s-1)\dots (s-m) ]^2  \  \Gamma(s+\nu) \left[  \Gamma(s)\right]^{k-2} y^{s-m-1} ds$$

$$= {1\over 2\pi i}  \int_{\gamma-i\infty}^{\gamma+i\infty} [(s)_m ]^2  \  \Gamma(s+m+\nu) \left[  \Gamma(s+m)\right]^{k-2} y^{s-1} ds$$

$$=  {y^{-m} \over 2\pi i}  \int_{\gamma-i\infty}^{\gamma+i\infty}   {\Gamma(s+\nu)\over \Gamma^2(s-m)}  \left[  \Gamma(s)\right]^{k} y^{s-1} ds.$$
Therefore we find from (4.15) 

$$ \int_0^\infty     Q_n^{\nu,\alpha,2} \left( y\right)   \Phi^{(2)}_{\nu,k,m} (y) \ y^{\alpha} dy  =  0,\ m= 0, 1,\dots, n-1,\   n \in \mathbb{N},\eqno(4.17)$$
where

$$ \Phi^{(2)}_{\nu,k,m} (y) \equiv  {1  \over 2\pi i}  \int_{\gamma-i\infty}^{\gamma+i\infty}   {\Gamma(s+m+\nu)\over \Gamma^2(s)}  \left[  \Gamma(s+m)\right]^{k} y^{-s} ds,\quad  k \ge 2.\eqno(4.18)$$
But it is easily seen from the properties of the Mellin transform \cite{YaL} and (2.1) that 

$$ \Phi^{(2)}_{\nu,k,m} (y) =  y^m D^m y^m D^m  y^m \left\{ \rho_{k-2} (y) \right\} ,\quad  k \ge 2.\eqno(4.19)$$
Now, recalling (4.4), we have 

$$y^m D^m y^m D^m y^m  \left\{ \rho_{k-2} (y) \right\} = y^m D^m y^m D^m y^m  \left\{  \int_0^\infty e^{-y u} \rho_{\nu,k-3} \left({1\over u}\right) {du\over u}\right\} .\eqno(4.20)$$
Hence, modifying the formula (4.14), we obtain

$$y^m D_y^m y^m D_y^m y^m \left\{ e^{-yu}\right\} = (-1)^m  D_u^m u^m D_u^m  u^m D_u^m \left\{ e^{-yu}\right\}.\eqno(4.21)$$
Therefore, integrating by parts, we get from (4.19), (4.20), (4.21)

$$  \Phi^{(2)}_{\nu,k,m} (y) =    \int_0^\infty   e^{-y u}   D_u^m u^m D_u^m  u^m D_u^m \left\{  \rho_{\nu,k-3} \left({1\over u}\right) {1\over u} \right\} du.\eqno(4.22)$$
Moreover, in a similar manner as above we derive

$$D_u^m u^m D_u^m  u^m D_u^m \left\{  \rho_{\nu,k-3} \left({1\over u}\right) {1\over u} \right\}  = D_u^m u^m D_u^m  u^m D_u^m  \left\{ {1\over 2\pi i}  \int_{\gamma-i\infty}^{\gamma+i\infty}  \Gamma(s+\nu) \left[  \Gamma(s)\right]^{k-3} u^{s-1} ds\right\} $$
 
$$= {1\over 2\pi i}  \int_{\gamma-i\infty}^{\gamma+i\infty} [(s-1)\dots (s-m) ]^3  \  \Gamma(s+\nu) \left[  \Gamma(s)\right]^{k-3} u^{s-m-1} ds$$

$$= {1\over 2\pi i}  \int_{\gamma-i\infty}^{\gamma+i\infty} [(s)_m ]^3  \  \Gamma(s+m+\nu) \left[  \Gamma(s+m)\right]^{k-3} u^{s-1} ds$$

$$= {u^{-m} \over 2\pi i}  \int_{\gamma-i\infty}^{\gamma+i\infty}  {\Gamma(s+\nu)\over \Gamma^3 (s-m) }\left[  \Gamma(s)\right]^{k} u^{s-1} ds.\eqno(4.23)$$
So, substituting the right-hand side of the latter equality in (4.23) into (4.22) and the obtained expression into (4.17), we find after the interchange of the order of integration and simple change of variables the following orthogonality conditions 

$$ \int_0^\infty     Q_n^{\nu,\alpha,3} \left( u\right)   \Phi^{(3)}_{\nu,k,m} (u) \ u^{\alpha} du  =  0,\ m= 0, 1,\dots, n-1,\   n \in \mathbb{N},\eqno(4.24)$$
where

$$ \Phi^{(3)}_{\nu,k,m} (u) \equiv  {1  \over 2\pi i}  \int_{\gamma-i\infty}^{\gamma+i\infty}   {\Gamma(s+m+\nu)\over \Gamma^3(s)}  \left[  \Gamma(s+m)\right]^{k} u^{-s} ds,\quad  k \ge 3,\eqno(4.25)$$
and 

$$ Q_n^{\nu,\alpha,3} \left( u\right) = {1\over \Gamma(1+\alpha)} \sum_{j=0}^n a_{n,j} [\Gamma (1+\alpha+j) ]^3  \ u^j.\eqno(4.26)$$
Continuing this process by virtue of  the same technique, involving the Mellin and Laplace transforms and the Mellin-Barnes integrals,  after  the  $k$-th step we end up  with the equalities  

$$ \int_0^\infty     Q_n^{\nu,\alpha,k} \left( x\right)   \Phi^{(k)}_{\nu,k,m} (x) \ x^{\alpha} dx  =  0,\ m= 0, 1,\dots, n-1,\   n \in \mathbb{N},\eqno(4.27)$$
where

$$ \Phi^{(k)}_{\nu,k,m} (x) =  {1  \over 2\pi i}  \int_{\gamma-i\infty}^{\gamma+i\infty}    [(s)_m)]^k \ \Gamma(s+m+\nu) x^{-s} ds,\eqno(4.28)$$
and 

$$ Q_n^{\nu,\alpha,k} \left( x\right) = {1\over \Gamma(1+\alpha)} \sum_{j=0}^n a_{n,j} [\Gamma (1+\alpha+j) ]^k  \ x^j.\eqno(4.29)$$
On the other hand,

$$ \Phi^{(k)}_{\nu,k,m} (x) =  {1  \over 2\pi i}  \int_{\gamma-i\infty}^{\gamma+i\infty}    [(s)_m)]^k \ \Gamma(s+m+\nu) x^{-s} ds
= (-1)^{km}  \left\{ x^m D^m \right\}^k \left( x^{\nu+m}  e^{-x} \right)$$

$$ =  (-1)^{km} m! \left\{ x^m D^m \right\}^{k-1}  \left( x^{\nu+m}  e^{-x} L_m^\nu(x) \right). $$
Consequently, the orthogonality (4.11) is equivalent to the following conditions

$$ \int_0^\infty     Q_n^{\nu,\alpha,k} \left( x\right)    \ x^{\alpha}  \left\{ x^m D^m \right\}^k \left( x^{\nu+m}  e^{-x} \right)  dx  =  0,\ m= 0, 1,\dots, n-1,\   n \in \mathbb{N},\ k \in \mathbb{N}_0.\eqno(4.30)$$
Moreover, we see that  $ \left\{ x^m D^m \right\}^k \left( x^{\nu+m}  e^{-x} \right) = x^{\nu}  e^{-x} p_{m(k+1)}(x),$ where $p_{m(k+1)}$ is a polynomial of degree $m(k+1)$, whose coefficients can be calculated explicitly via properties of  the Pochhammer symbol and the associated Laguerre polynomials.  Thus it can be reduced to the orthogonality with respect to the measure $x^{\nu+\alpha} e^{-x} dx$ and ideas of the previous section can be applied. We leave all details to  the interested reader.  Besides, further developments, an analog of Lemma 1 and relations with the multiple orthogonal polynomial ensemble  from \cite{Arno} will be a promising  investigation.

\bibliographystyle{amsplain}

\begin{thebibliography}{10}

\bibitem{Bateman} A. Erd\'elyi,W. Magnus, F. Oberhettinger, and F.G. Tricomi, Higher Transcendental Functions,Vols. I and II, McGraw-Hill, NewYork, London, Toronto, 1953.
\bibitem{Chi} T.S. Chihara, An Introduction to Orthogonal Polynomials, Gordon and Breach, New York, London, 1978. 
\bibitem{Cous}  E. Coussement, W. Van Assche, Some properties of multiple orthogonal polynomials associated with Macdonald functions,     J. Comput. Appl. Math. 133 (2001),  253-261.
\bibitem{Dattoli}  G. Dattoli, P.E. Ricci, I. Khomasuridze,  Operational methods, special polinomial and functions and solution of partial differential equations,   Integral Transforms and  Special  Functions,  15 (2004), N 4,  309- 321.
\bibitem{DitkinPrudnikov} V. A. Ditkin,  A.P. Prudnikov, Integral transforms. Mathematical analysis, 1966,  7-82. Akad. Nauk SSSR Inst. Nauchn. Informacii, Moscow, 1967 (in Russian).
\bibitem{Ismail}  M.E.H. Ismail,  Bessel functions and the infinite divisibility of the student $t$-distribution, Ann. Prob. 5 (1977), 582-585. 
\bibitem{Arno}  A.B.J. Kuijlaars,  L. Zhang,    Singular values of products of Ginibre random matrices, multiple orthogonal polynomials and hard eddge scaling limits,      Commun. Math. Phys. 332  (2014),  759-781.
\bibitem{Filipa}  A.F. Loureiro, S.Yakubovich,  Central factorials under the Kontorovich-Lebedev transform of polynomials,   Integral Transforms and  Special  Functions,  24 (2013), N 3,  217-238.
\bibitem{PrudnikovProblem} A.P. Prudnikov, Orthogonal polynomials with ultra-exponential weight functions, in: W. Van Assche (Ed.), Open Problems, J. Comput. Appl. Math. 48 (1993),  239-241.
\bibitem{PrudnikovMarichev} A. P. Prudnikov, Yu. A. Brychkov and O.I. Marichev, Vol. I: Elementary Functions, Gordon and Breach, New York, London, 1986; Vol. II: Special Functions,  Gordon and Breach, New York, London, 1986;   Vol. III: More Special Functions, Gordon and Breach, New York, London, 1990. 
\bibitem{Riordan} J. Riordan, {Combinatorial identities}, Wiley, New Yok, 1968.
\bibitem{Tit} E.C. Titchmarsh,  An Introduction to the Theory of Fourier Integrals, Clarendon Press, Oxford, 1937.
\bibitem{AsscheYakubov2000} W. Van Assche and S. Yakubovich,  Multiple orthogonal polynomials associated
with Macdonald functions, Integral Transforms and  Special  Functions,  9 (2000), N 3, 229- 244.
\bibitem{Viskov}  O.V. Viskov, H.M. Srivastava,  New approaches to certain identities involving differential operators, J. Math. Anal. Appl. 186 (1994), 1-10. 
\bibitem{YaL} S. Yakubovich and Yu. Luchko, The Hypergeometric Approach to Integral Transforms and Convolutions, Kluwer
Academic Publishers, Mathematics and Applications. Vol.287, 1994. 
\end{thebibliography}

\end{document}